\documentclass[11pt, a4paper]{article}
\usepackage[textwidth = 7in, textheight=25cm,nohead]{geometry}
\usepackage{ulem}
\usepackage{arydshln}
\usepackage{graphicx, amssymb, amsmath, amsthm, amstext, booktabs, float, multirow, rotating,natbib,bm}
\usepackage{caption}
\usepackage{subcaption}
\usepackage{stackengine}
\usepackage{amssymb}
\usepackage{mathtools,bm}
\usepackage[mathscr]{euscript}
\usepackage{upgreek}

\usepackage{bbding}
\usepackage{pifont}

\usepackage{xparse}

\usepackage{mathdots}

\let\emptyset\varnothing

\usepackage{tikz}
\usepackage{setspace}
\usepackage{enumerate}
\usepackage{multirow}
\usepackage{graphicx,nicefrac}
\newcommand\bsfrac[2]{%
\scalebox{-1}[1]{\nicefrac{\scalebox{-1}[1]{$#1$}}{\scalebox{-1}[1]{$#2$}}}%
}

\numberwithin{equation}{section} 
\numberwithin{figure}{section} 
\numberwithin{table}{section} 
\newtheorem{definition}{Definition}[section] 
\newtheorem{lemma}{Lemma}[section] 

\newtheorem{theorem}{Theorem}[section] 
 
\newtheorem{proposition}{Proposition}[section] 
\newtheorem{corollary}{Corollary}[section]

\theoremstyle{remark} 
\newtheorem{remark}{Remark}[section] 

\theoremstyle{definition} 
\newtheorem{fact}{Fact}[section] 

\theoremstyle{definition} 
\newtheorem{example}{Example}[section]

\bibpunct{(}{)}{,}{a}{,}{,}

\makeatletter
\newcommand{\doublewidetilde}[1]{{%
\mathpalette\double@widetilde{#1}%
}}
\newcommand{\double@widetilde}[2]{%
\sbox\z@{$\m@th#1\widetilde{#2}$}%
\ht\z@=.9\ht\z@
\widetilde{\box\z@}%
}
\makeatother

 \let\mathscr\relax% just so we can load this and rsfs
\usepackage[scr]{rsfso}
\newcommand{\powerset}{\raisebox{.15\baselineskip}{\Large\ensuremath{\wp}}}

\def\ppn{\vskip 6pt \noindent }
\def\R{{\mathbb{R}}}
\def\N{{\mathbb{N}}}

\def\P{{\mathbb{P}}}
\def\E{{\mathbb{E}}}
\newcommand{{\Xs}}{{\cal X}}
\newcommand{{\Xss}}{{\mathfrak{X}}}
\newcommand{{\Ys}}{{\cal Y}}
\newcommand{{\Yss}}{{\mathfrak{Y}}}
\newcommand{{\Ls}}{{\cal L}}
\newcommand{{\Ss}}{{\cal S}}
\newcommand{{\Ms}}{{\cal M}}
\newcommand{{\Gs}}{{\cal G}}
\newcommand{{\Hs}}{{\cal H}}
\newcommand{{\Ns}}{{\cal N}}
\newcommand{{\Is}}{{\cal I}}
\newcommand{{\Vs}}{{\cal V}}
\newcommand{{\Ds}}{{\cal D}}
\newcommand{{\As}}{{\cal A}}
\newcommand{{\Bs}}{{\cal B}}
\newcommand{{\Cs}}{{\cal C}}
\newcommand{{\Rs}}{{\cal R}}
\newcommand{{\Es}}{{\cal E}}
\newcommand{{\Fs}}{{\cal F}}
\newcommand{{\Us}}{{\cal U}}
\newcommand{{\Ps}}{{\cal P}}
\newcommand{{\ttheta}}{{\bm{\theta}}}
\newcommand{{\Ttheta}}{{\bm{\Theta}}}
\newcommand{{\Oomega}}{{\bm{\Omega}}}
\newcommand{{\oomega}}{{\bm{\omega}}}
\newcommand{{\mmu}}{{\bm{\mu}}}
\newcommand{{\ggamma}}{{\bm{\gamma}}}
\newcommand{{\Ssigma}}{{\bm{\Sigma}}}
\newcommand{{\Sss}}{{\bm{\Ss}}}
\newcommand{{\pp}}{{\mathbf p}}
\newcommand{{\ww}}{{\mathbf w}}
\newcommand{{\mm}}{{\mathbf m}}
\newcommand{{\bb}}{{\mathbf b}}
\newcommand{{\uu}}{{\mathbf u}}
\newcommand{{\ppi}}{{\bm{\pi}}}
\newcommand{{\phhi}}{{\bm{\phi}}}
\newcommand{{\pssi}}{{\bm{\psi}}}
\newcommand{{\XX}}{{\mathbf X}}
\newcommand{{\UU}}{{\mathbf U}}
\newcommand{{\BB}}{{\mathbf B}}
\newcommand{{\KK}}{{\mathbf K}}
\newcommand{{\HH}}{{\mathbf H}}
\newcommand{{\II}}{{\mathbf I}}
\newcommand{{\PP}}{{\mathbf P}}
\newcommand{{\yy}}{{\mathbf y}}
\newcommand{{\ee}}{{\mathbf e}}
\newcommand{{\ab}}{{\mathbf a}}

\newcommand{{\dd}}{{\mathbf d}}
\newcommand{{\zero}}{{\mathbf 0}}
\newcommand{{\uno}}{{\mathbf 1}}

\def\dep{\mathfrak{D}}

\newcommand{{\DDelta}}{{\bm \Delta}}
\newcommand{{\piXY}}{{\pi_{{\scriptscriptstyle XY}}}}
\newcommand{{\FXY}}{{F_{{\scriptscriptstyle XY}}}}
\newcommand{{\fXY}}{{f_{{\scriptscriptstyle XY}}}}
\newcommand{{\piXYstar}}{{\pi^*_{{\scriptscriptstyle XY}}}}
\newcommand{{\piXYstarstar}}{{\pi^{**}_{{\scriptscriptstyle XY}}}}
\newcommand{{\piotstar}}{{\pi^*_{{\scriptscriptstyle 12}}}}
\newcommand{{\piX}}{{\pi_{{\scriptscriptstyle X}}}}
\newcommand{{\FX}}{{F_{{\scriptscriptstyle X}}}}
\newcommand{{\fX}}{{f_{{\scriptscriptstyle X}}}}
\newcommand{{\pio}}{{\pi_{{\scriptscriptstyle 1}}}}
\newcommand{{\piXstar}}{{\pi^*_{{\scriptscriptstyle X}}}}
\newcommand{{\piostar}}{{\pi^*_{{\scriptscriptstyle 1}}}}
\newcommand{{\piY}}{{\pi_{{\scriptscriptstyle Y}}}}
\newcommand{{\FY}}{{F_{{\scriptscriptstyle Y}}}}
\newcommand{{\fY}}{{f_{{\scriptscriptstyle Y}}}}
\newcommand{{\pit}}{{\pi_{{\scriptscriptstyle 2}}}}
\newcommand{{\piYstar}}{{\pi^*_{{\scriptscriptstyle Y}}}}
\newcommand{{\pitstar}}{{\pi^*_{{\scriptscriptstyle 2}}}}
\newcommand{{\phiXY}}{{\varphi_{{\scriptscriptstyle XY}}}}
\newcommand{{\psiXY}}{{\psi_{{\scriptscriptstyle XY}}}}
\newcommand{{\phiXYstar}}{{\varphi^*_{{\scriptscriptstyle XY}}}}
\newcommand{{\phiX}}{{\varphi_{{\scriptscriptstyle X}}}}
\newcommand{{\phiXstar}}{{\varphi^*_{{\scriptscriptstyle X}}}}
\newcommand{{\phiYstar}}{{\varphi^*_{{\scriptscriptstyle Y}}}}
\newcommand{{\phiY}}{{\varphi_{{\scriptscriptstyle Y}}}}
\newcommand{{\phiYgX}}{{\varphi_{{\scriptscriptstyle Y|X}}}}
\newcommand{{\phiXgY}}{{\varphi_{{\scriptscriptstyle X|Y}}}}
\newcommand{{\phiXgYstar}}{{\varphi^*_{{\scriptscriptstyle X|Y}}}}
\newcommand{{\phiYgXstar}}{{\varphi^*_{{\scriptscriptstyle Y|X}}}}

\newcommand\indep{\protect\mathpalette{\protect\independenT}{\perp}}
\def\independenT#1#2{\mathrel{\rlap{$#1#2$}\mkern2mu{#1#2}}}

\newcommand*\diff{\mathop{}\!\mathrm{d}} %% places a "d" in dx

\newcommand{{\toL}}{{\overset{\mathcal{L}}{\longrightarrow}\ }}

\newcommand{{\MC}}{{\,  *_{\text{\scalebox{0.65}{$\Ms$}}}\,  }}

\newcommand{{\dou}}{$\leadsto$\ }

\DeclareMathOperator{\Ran}{Ran}

\DeclareMathOperator{\atanh}{atanh}
\DeclareMathOperator{\asin}{asin}

\newcommand{\indic}[1]{
\hbox{${\it 1}\hskip -4.5pt I_{\{ #1 \}}$}
}

\begin{document}

\setlength{\belowdisplayskip}{5pt} \setlength{\belowdisplayshortskip}{3pt}
\setlength{\abovedisplayskip}{5pt} \setlength{\abovedisplayshortskip}{0pt}

\title{Towards a universal representation of statistical dependence}
\author{\sc{Gery Geenens}\thanks{Corresponding author: {\tt ggeenens@unsw.edu.au}, School of Mathematics and Statistics, UNSW Sydney, Australia, tel +61 2 938 57032, fax +61 2 9385 7123 }\\School of Mathematics and Statistics,\\ UNSW Sydney, Australia 
}
\date{\today}
\maketitle
\thispagestyle{empty} 

\begin{abstract} \noindent Dependence is undoubtedly a central concept in statistics. Though, it proves difficult to locate in the literature a formal definition which goes beyond the self-evident `dependence = non-independence'. This absence has allowed the term `dependence' and its declination to be used vaguely and indiscriminately for qualifying a variety of disparate notions, leading to numerous incongruities. For example, the classical Pearson's, Spearman's or Kendall's correlations are widely regarded as `dependence measures' of major interest, in spite of returning 0 in some cases of deterministic relationships between the variables at play -- evidently not measuring dependence at all. Arguing that research on such a fundamental topic would benefit from a slightly more rigid framework, this paper suggests a general definition of the dependence between two random variables defined on the same probability space. Natural enough for aligning with intuition, that definition is still sufficiently precise for allowing unequivocal identification of a `universal' representation of the dependence structure of any bivariate distribution. Links between this representation and familiar concepts are highlighted, and ultimately, the idea of a dependence measure based on that universal representation is explored and shown to satisfy R\'enyi's postulates.
\end{abstract}

%\ppn {\bf Keywords:} %Copula; discrete random vector; bivariate Bernoulli distribution; Yule's colligation coefficient; Iterative Proportional Fitting.

%\ppn {\bf AMS Classification:} %62H05, 62H17, 62H20.

\section{Introduction}\label{sec:intro}

Scientific progress relies on identifying if, how, why, and to which extent, several factors influence each other. Across all fields, studies looking for establishing an effect of an input on an output must statistically demonstrate evidence of association\,/\,correlation\,/\,relationship\,/\,causality between input and output; i.e., to establish some sort of {dependence} between the two. Hence statistics, the art of turning empirical evidence ({data}) into information, has always kept the concept of dependence at its core.

\ppn In spite of this, it may prove difficult to delineate `dependence' precisely. It is noteworthy that, notwithstanding the frequent use of the term, it is rarely defined formally in the literature -- not even in major references entirely dedicated to the topic such as \cite{DrouetMari01} or \cite{Joe2015}. Neither the `{Cambridge Dictionary of Statistics}' \citep{Everitt10} nor the `{Dictionary of Statistics}' \citep{Upton14} include an entry for dependence. The `{Dictionary of Statistical Terms}' \citep{Kendall71,Dodge06} does, but gives the seemingly circular `{\it Dependence: quantities are dependent when they are not independent}'. 

\ppn Admittedly, {\it independence} has an unambiguous mathematical meaning, as when the joint distribution of two variables of interest can be suitably factorised; see (\ref{eqn:indep2}) below. Writing this in terms of conditional distributions allows an interpretation along the line of `variation in one variable does not disturb the stochastic behaviour of the other'. Thus, according to the above (pseudo-)definition, two variables are dependent as soon as {there is} such (loosely defined, not necessarily causal) influence of one on the other, and dependence amounts to the {\it existence} of the said influence. In effect such `{\it non-independence}' interpretation of dependence makes it a binary concept -- two variables do or do not satisfy the property of independence, so they are dependent or not without any room for nuance. If such a binary-in-nature view may be suitable for the {\it independence testing} problem, it quickly falls short when it comes to other questions. In particular, quantifying dependence appears essential in many situations, as testified by the abundance of {dependence measures} found in the literature -- see \cite{Tjostheim22} for a recent review. Though, measuring a trait is incompatible with an alleged binary nature. If dependence is to be quantifiable, then clearly the non-independence definition appears inadequate -- and this, without any obvious substitute.

\ppn This absence has permitted some looseness when addressing dependence-related questions, inducing confusion between disparate notions, and causing incongruities. The most obvious of those may be that the most popular `dependence measures', namely Pearson's correlation ($\rho$), Spearman's rho ($\rho_S$) and Kendall's tau ($\tau$), do {\it not} measure dependence at all -- for example, they may be 0 between two variables deterministically bound,\footnote{See the textbook example $X \sim \Ns(0,1)$, $Y=X^2$ $\Rightarrow$ $\rho(X,Y) = \rho_S(X,Y) = \tau(X,Y) = 0 $.} i.e., the strongest influence that one can imagine. In fact, the elusiveness of a quantifiable version of dependence has allowed for subjective interpretation to resolve what was to be measured, and how. This explains the profusion of diverse, more or less pertinent, `dependence measures' (more~on~this~in~Section~\ref{sec:depmes}).

\ppn Arguing that further progress in an area as important as that of statistical dependence would benefit from a slightly more rigid framework, this paper attempts to provide a clearer understanding of the concept, starting from first principles (Section \ref{sec:heuristic}). We then propose a general definition of the dependence between two variables defined on the same probability space, which leads, implication after implication (through Section \ref{sec:dep}), to a unequivocal representation of the said concept (Section \ref{sec:universal0}). This representation is `universal' in the sense that it is valid for any random vector, regardless of its nature (discrete, continuous, mixed or hybrid) -- we only require a joint distribution dominated by a reference product measure. Section \ref{sec:concord} illustrates the distinction between dependence and other related notions by examining closely the important concept of concordance. Section \ref{sec:depmes} suggests a general measure of dependence articulated around the obtained universal representation, and Section \ref{sec:ccl} concludes. 

\section{Heuristic} \label{sec:heuristic}

Let $(X,Y)$ be a bivariate random vector with unknown joint distribution $\piXY$. Suppose that an oracle gives us full knowledge of the two univariate marginal distributions $\piX$ and $\piY$. Is it enough to reconstruct $\piXY$ entirely? The answer is no: this information allows us to conclude that $\piXY$ belongs to the {Fr\'echet class} $\Fs(\piX,\piY)$ -- the set of all bivariate distributions with margins $\piX$ and $\piY$ -- but surely such class contains more than one element. Then, what is the extra piece of information, call it $\dep$, that the oracle is still concealing, and which would allow unequivocal identification of $\piXY$ in $\Fs(\piX,\piY)$? 

\ppn Specifically, $\dep$ stands for what is {necessary} and {sufficient} for identifying $\piXY$ once $\piX$ and $\piY$ are known: `necessary' as it does not include any redundant information already contained in the margins, and `sufficient' as nothing relevant to $\piXY$ is possibly left out. This breaks down the bivariate distribution $\piXY$ into three constituents showing no overlap: $\piX$, describing the behaviour of $X$ alone; $\piY$, describing the behaviour of $Y$ alone; and $\dep$, describing `the rest'.

\ppn The marginal distributions $\piX$ and $\piY$ may be given various mathematical representations, such as cumulative distribution functions, probability mass/density functions or characteristic functions, the suitability of which depending on the nature of $X$ and $Y$ and the assumptions made {a priori} about them. E.g., if $\piX$ is supported on a finite set, then exhaustive enumeration of the relevant probabilities (`probability mass function') is possible, but not if $\piX$ has uncountable support; while if $\piX$ is assumed to be Gaussian, then the pair (mean, variance) provides an adequate description by itself. We anticipate that the third element $\dep$ is no different, with potentially several mathematical representations differing according to the nature of $(X,Y)$ -- this will be formally confirmed in the following sections. For now we may just understand $\dep$ informally as what remains of $\piXY$ when entirely stripped from its margins, i.e., the `glue' that keeps $\piX$ and $\piY$ together in $\piXY$.

\ppn Clearly, if $X$ and $Y$ happen to be independent, $\dep$ contains that information -- and only that, as it is enough for identifying $\piXY = \piX \times \piY$ in $\Fs(\piX,\piY)$. If $X$ and $Y$ are not independent, then $\dep$ should provide the complete recipe as to how to interlock $\piX$ and $\piY$ for giving rise to the specific $\piXY$ and making it different to any other distribution of $\Fs(\piX,\piY)$. In all cases, $\dep$ describes what happens `between' $\piX$ and $\piY$,  which it seems fair to call the dependence structure -- or just dependence -- of the distribution $\piXY$. Identifying a random vector to its distribution, we may also say that $\dep$ is the dependence of the vector $(X,Y)$. We, therefore, propose the following definition:

\begin{definition} \label{def:dep}
	The {\bf dependence} between two random variables $X \sim \piX$ and $Y \sim \piY$ defined on the same probability space is the information $\dep$ which is necessary and sufficient to unequivocally identify their joint distribution $\piXY$ in its Fr\'echet class $\Fs(\piX,\piY)$.
\end{definition}

\noindent At this stage, this definition is purposely qualitative and very general. Yet, through a sequence of logical implications, it will translate into a precise mathematical representation (Corollary \ref{cor:dep}) -- identifying such a representation from the intuitive Definition \ref{def:dep} is actually the main goal of this paper. As-is, it suggests the equivalence 
\begin{equation} \piXY \longleftrightarrow (\piX,\piY;\dep),\label{eqn:decomp} \end{equation}
where $\piXY \longrightarrow (\piX,\piY;\dep)$ means that any distribution $\piXY$ admits one and only one decomposition in terms of its three constituents $\piX$, $\piY$ and $\dep$; while $\piXY \longleftarrow (\piX,\piY;\dep)$ means that any triplet $(\piX,\piY;\dep)$ defines one and only one distribution $\piXY$ with the corresponding margins and dependence structure. The ($\longrightarrow$)-part is true because any distribution in a given Fr\'echet class must be unequivocally distinguished from other distributions in the same class by some features not perceptible from the marginals, establishing the existence of a third element $\dep$ specific to $\piXY$. In Section \ref{sec:dep} below, we prove that the ($\longleftarrow$)-part holds true, too, up to some natural conditions on the cardinality and shape of the supports of the distributions at play. E.g., as anticipated above, $\dep$ may be of different nature when $\piX,\piY$ are discrete and when they are continuous, hence a `continuous $\dep$' cannot be combined with discrete marginals (or vice-versa). But, dependence structures `of the right type' may be freely interchanged for cementing any given marginals $\piX, \piY$ and giving rise to various distributions of $\Fs(\piX,\piY)$.

\ppn For this to be possible, though, it is necessary that no specification of any of the three constituents alters or imposes restrictions on the others. In a sense, a distribution $\piXY$ should be a point represented by a triplet $(\piX,\piY,\dep)$ of {orthogonal coordinates} in a some sort of rectangular-prismatic Cartesian space ($X$-margin, $Y$-margin, dependence). This is comparable to the concept of {variation-independence} between parameters \cite[p.\ 26]{Barndorff78}: in a parametric model, two parameters $\theta$ and $\theta'$ are variation-independent if setting the value of one does not impose any new logical constraints on the range of the other. If the decomposition $(\piX,\piY;\dep)$ is regarded as a parametrisation of $\piXY$, then the above discussion implies that $\piX$, $\piY$ and $\dep$ should be variation-independent. In particular, if any specific feature of $\piXY$ happens not to be variation-independent of $\piX$ and $\piY$, then it cannot enter any valid description of $\dep$ as such since it would automatically carry redundant information about the margins, in violation of Definition \ref{def:dep}. Here three simple examples are used to illustrate the point.

\begin{example} \label{ex:Gauss} {\bf Bivariate Gaussian.} Let $(X,Y)$ be a bivariate Gaussian vector. Then $\piXY$ is entirely described by its mean vector and its variance-covariance matrix, that is, 
	\begin{equation} \mmu  = \begin{pmatrix} \mu_{\scriptscriptstyle X} \\ \mu_{\scriptscriptstyle Y} \end{pmatrix}, \qquad \Ssigma = \begin{pmatrix} \sigma^2_{\scriptscriptstyle X} & \sigma_{\scriptscriptstyle XY} \\ \sigma_{\scriptscriptstyle XY} & \sigma^2_{\scriptscriptstyle Y}\end{pmatrix}  = \begin{pmatrix} \sigma^2_{\scriptscriptstyle X} & \rho \sigma_{\scriptscriptstyle X} \sigma_{\scriptscriptstyle Y} \\ \rho \sigma_{\scriptscriptstyle X} \sigma_{\scriptscriptstyle Y} & \sigma^2_{\scriptscriptstyle Y}\end{pmatrix}. \label{eqn:paramGauss} \end{equation}	
	Once $\piX$ and $\piY$ have been specified -- i.e., $X \sim \Ns(\mu_{\scriptscriptstyle X},\sigma^2_{\scriptscriptstyle X})$ and $Y \sim \Ns(\mu_{\scriptscriptstyle Y},\sigma^2_{\scriptscriptstyle Y})$ -- the only free parameter left in (\ref{eqn:paramGauss}) is $\rho = \frac{\sigma_{\scriptscriptstyle XY}}{\sigma_{\scriptscriptstyle X} \sigma_{\scriptscriptstyle Y}}$ (Pearson's correlation). As such, $\rho$ alone allows unequivocal identification of $\piXY$ {\it among the bivariate Gaussian distributions} of $\Fs(\Ns(\mu_{\scriptscriptstyle X},\sigma^2_{\scriptscriptstyle X}),\Ns(\mu_{\scriptscriptstyle Y},\sigma^2_{\scriptscriptstyle Y}))$ but does not contain any information about the marginal parameters. Therefore, $\rho$ 
	appropriately captures {\it the whole dependence structure} of $(X,Y)$, as per Definition \ref{def:dep} -- and is not just a quantification of the strength of association. Irrespective of $(\mu_{\scriptscriptstyle X},\sigma_{\scriptscriptstyle X})$ and $(\mu_{\scriptscriptstyle Y},\sigma_{\scriptscriptstyle Y})$, for any $\rho \in [-1,1]$ there exists one and only one bivariate Gaussian distribution with those margins and correlation $\rho$; whence $(\mu_{\scriptscriptstyle X},\sigma_{\scriptscriptstyle X})$, $(\mu_{\scriptscriptstyle Y},\sigma_{\scriptscriptstyle Y})$ and $\rho$ are indeed variation-independent. Obviously, any one-to-one function of $\rho$ would allow identification of $\piXY$ as well: it is the case of $z=\atanh \rho$ (Fisher's $z$-transform), Kendall's $\tau = \frac{2}{\pi} \asin \rho$, and Spearman's $\rho_S = \frac{6}{\pi} \asin \frac{\rho}{2}$ \citep{Kruskal58}, which therefore qualify all the same to represent the dependence here. By contrast, the covariance $\sigma_{\scriptscriptstyle XY}$ does not, as it contains redundant marginal information: $\sigma_{\scriptscriptstyle X}$, $\sigma_{\scriptscriptstyle Y}$ and $\sigma_{\scriptscriptstyle XY}$ are not variation-independent, as $|\sigma_{\scriptscriptstyle XY}| \leq \sigma_{\scriptscriptstyle X} \sigma_{\scriptscriptstyle Y}$. %Concretely, if an oracle lets us know that $\sigma_{\scriptscriptstyle XY}=10$ in a bivariate Gaussian distribution, this does not tell us anything about the dependence as long as we do not consider this value relatively to the marginal parameters $\sigma_{\scriptscriptstyle X}$ and $\sigma_{\scriptscriptstyle Y}$.
	
	\ppn Naturally the parity `$\dep$ $\sim$ Pearson's correlation $\rho$' follows from the rigid specification of the Gaussian model, but breaks as soon as we leave it. This explains why $\rho = 0 \iff X \indep Y$ (independence) if $(X,Y)$ is Gaussian, but not in general. In fact, the value of $\rho$ does not unequivocally single out a distribution in the whole of $\Fs(\Ns(\mu_{\scriptscriptstyle X},\sigma^2_{\scriptscriptstyle X}),\Ns(\mu_{\scriptscriptstyle Y},\sigma^2_{\scriptscriptstyle Y}))$, as this class contains non-Gaussian elements with identical correlation but different dependence structures. E.g., there exist distributions with Gaussian marginals and null correlation, which are not the bivariate Gaussian distribution with independent components.\qed
\end{example}

\begin{example} \label{ex:bivBern} {\bf Bivariate Bernoulli.} Let $(X,Y)$ be a bivariate Bernoulli vector; that is, for $p_{1\bullet},p_{\bullet 1} \in (0,1)$, $X \sim \text{Bern}(p_{1\bullet})$, $Y \sim \text{Bern}(p_{\bullet 1})$ and the joint distribution $\piXY$ is described by a $(2 \times 2)$-table like
	\begin{equation} \begin{array}{c l ||c c c c | c}
			& \bsfrac{\ Y}{X\ } & & 0 & 1 & \\
			\hline\hline
			& 0 & & p_{00} & p_{01}& & p_{0\bullet} \\
			& 1 & & p_{10}  & p_{11}& & p_{1\bullet}  \\
			\hline 
			& & & p_{\bullet 0} & p_{\bullet 1}& & 1 \\
		\end{array}, \label{eqn:tab22} \end{equation}
	where\footnote{Here we consider only full tables (all positive entries) for simplicity. Presence of structural zeros is treated in full generality in the following sections.}  $p_{xy} = \P(X=x,Y=y) >0$, $p_{x \bullet} = p_{x0}+p_{x1}$ and $p_{\bullet y} = p_{0y} + p_{1y}$ ($x,y \in \{0,1\}$). As $\sum p_{xy} = 1$, there are initially 3 free parameters for $\piXY$. Once we specify the margins $\piX$ and $\piY$, that is, once we fix the two values $p_{1\bullet}$ and $p_{\bullet 1}$, only one free parameter is left -- the one which would identify $\piXY$ among all the distributions of $\Fs(\text{Bern}(p_{1\bullet}),\text{Bern}(p_{\bullet 1 }))$. Thus, as in the Gaussian framework (Example \ref{ex:Gauss}), the dependence $\dep$ is fully represented by one single number. However, that number cannot be Pearson's correlation which is here:
	\[ \rho = \frac{p_{11}-p_{1\bullet}p_{\bullet 1} }{\sqrt{p_{1\bullet}(1-p_{1\bullet})p_{\bullet 1}(1-p_{\bullet 1})}}.\]
	Indeed, $p_{11} \in \left( \max(0,p_{1\bullet} + p_{\bullet 1}-1),\min(p_{1\bullet},p_{\bullet 1})\right)$ -- the so-called `Fr\'echet bounds' \citep{Frechet51,Hoeffding40} in $\Fs(\text{Bern}(p_{1\bullet}),\text{Bern}(p_{\bullet 1}))$ -- which in turn imposes
	\[-\left(\min\{\frac{p_{1\bullet}p_{\bullet 1}}{(1-p_{1\bullet})(1-p_{\bullet 1})},\frac{(1-p_{1\bullet})(1-p_{\bullet 1})}{p_{1\bullet}p_{\bullet 1}}\}\right)^{1/2}< \rho < \left(\min\{\frac{p_{1\bullet}(1-p_{\bullet 1})}{(1-p_{1\bullet})p_{\bullet 1}},\frac{p_{\bullet 1}(1-p_{ 1 \bullet})}{(1-p_{\bullet 1})p_{ 1 \bullet}}\}\right)^{1/2}.\]
	So, not all combinations $(p_{1 \bullet},p_{\bullet 1},\rho) \in (0,1) \times (0,1) \times (-1,1)$ produce a valid bivariate Bernoulli distribution, and $\rho$ is not variation-independent of the marginal parameters $p_{1\bullet}$ and $p_{\bullet 1}$. The same holds for any quantity related to the $\chi^2$-statistic such as the contingency coefficient $\Phi^2$, or any one of the `association factors'/`Pearson ratios' \citep{Good56,Goodman96}:
	\begin{equation} \psi_{xy} \doteq \frac{p_{xy}}{p_{x \bullet}p_{\bullet y}} \in  \left(\max(0,\frac{1}{p_{x\bullet}} + \frac{1}{p_{\bullet y}} - \frac{1}{p_{x\bullet} p_{\bullet y}}) \ ,\ \frac{1}{\max(p_{x\bullet}, p_{\bullet y})}\right) \qquad (x,y) \in \{0,1\} \times \{0,1\}. \label{eqn:psi} \end{equation}
	Likewise, the above bounds constrain Spearman's $\rho_S$ and Kendall's $\tau$, both equal to $p_{11} - p_{1\bullet} p_{\bullet 1}$ in this case \citep[Example 7]{Genest07}, to be such that
	\[ \max\big(-p_{1\bullet}p_{\bullet 1},-(1-p_{1\bullet})(1-p_{\bullet 1})\big) \leq \rho_S, \tau \leq \min(p_{1\bullet},p_{\bullet 1}) \big(1-\max(p_{1\bullet},p_{\bullet 1})\big),\]
	thus $\rho_S$ and  $\tau$ are not variation-independent to $p_{1\bullet}$ and $p_{\bullet 1}$.\footnote{This remains true for any `tie-corrected' version such as $\tau_b$ \citep[Chapter 3]{Kendall55} or $\tau_c$ \citep{Stuart53}.} Consequently, none of these can be thought of as appropriately describing $\dep$ in a distribution such as (\ref{eqn:tab22}). 
	
	\ppn In fact, results of \cite{Edwards63} and \citet[Theorem 6.3]{Rudas18} establish that any parameter allowing identification of the distribution (\ref{eqn:tab22}) in $\Fs(\text{Bern}(p_{1\bullet}),\text{Bern}(p_{\bullet 1 }))$, while being variation-independent of $(p_{1 \bullet},p_{\bullet 1})$, must be a one-to-one function of the odds-ratio%\footnote{Here we allow $\omega = \infty$, in cases $p_{10}=0$ or $p_{01}=0$.}
	\begin{equation*} \omega = \frac{p_{00}p_{11}}{p_{10}p_{01}} \in (0,\infty). \label{eqn:OR} \end{equation*}
	Indeed any triplet $(p_{1\bullet},p_{\bullet 1},\omega) \in (0,1) \times (0,1) \times (0,\infty)$ defines one and only one table (\ref{eqn:tab22}) \cite[equation (2*)]{Mosteller68a}. Any one-to-one function of $\omega$ is thus a valid representation of $\dep$ as well; this includes $\log \omega$, \cite{Yule12}'s `association' and `colligation' coefficients $Q = (\omega -1)/(\omega+1)$ and $\Upsilon = (\sqrt{\omega}-1)/(\sqrt{\omega}+1)$. 	\qed
\end{example}

\begin{example} \label{ex:bivCop} {\bf Bivariate copula model.} Let the cumulative distribution function (cdf) of $\piXY$ be $\FXY$, and the cdf's of $\piX$ and $\piY$ be $\FX$ and $\FY$, respectively. It follows from \citet[{\it Th\'eor\`eme} 2]{Sklar59} that there exists a {\it copula} $C$ -- that is, a bivariate cdf supported on $[0,1]^2$ with standard uniform margins -- such that 
	\begin{equation} \FXY(x,y) = C(\FX(x),\FY(y)), \qquad \forall (x,y). \label{eqn:Sklar} \end{equation}
	This establishes such $C$ as a natural candidate for representing the dependence $\dep$ in any bivariate distribution, as it appears to capture how to reconstruct $\FXY$ from $\FX$ and $\FY$. If $X$ and $Y$ are both continuous variables, $C$ is unique and can easily be seen to be the distribution of the vector $(\FX(X),\FY(Y))$ \citep[Theorem 3$(iii)$]{Schweizer81}. Therefore, by the virtue of the Probability Integral Transform result ($\FX(X),\FY(Y) \sim \Us_{[0,1]}$), it is `marginal-distribution-free',\footnote{Where `distribution-free' is understood in the sense of \cite{Kendall53}: {\it free of the parent distribution}. Thus, more specifically here, `marginal-distribution-free'  means {\it free of the marginal distributions of the parent distribution $\piXY$}.} that is, `margin-free'. Thus, $C$ provides the margin-free information which allows unequivocal identification of $\piXY$ in $\Fs(\piX,\piY)$, which is $\dep$ by Definition \ref{def:dep}. This is consistent with Example \ref{ex:Gauss}, as the copula of a bivariate Gaussian distribution is by definition the Gaussian copula \citep[Section 4.3.1]{Joe2015}, whose only free parameter is the correlation $\rho$. By contrast, the dependence of the non-Gaussian distributions of $\Fs(\Ns(\mu_{\scriptscriptstyle X},\sigma^2_{\scriptscriptstyle X}),\Ns(\mu_{\scriptscriptstyle Y},\sigma^2_{\scriptscriptstyle Y}))$ is characterised by other copulas, not driven by $\rho$.
	
	\ppn In non-continuous cases, though, the copula is not unique, meaning that any $C$ satisfying (\ref{eqn:Sklar}) automatically contains irrelevant information not related to $\dep$. Evidently, what (\ref{eqn:Sklar}) may uniquely identify is only the restriction of $C$ to\footnote{$\Ran$ denotes the range: $\Ran F = \{t \in [0,1]: \exists x \in \overline{\R} \text{ s.t. } F(x) = t\}$.} $\Ran \FX \times \Ran \FY$ \citep[{\it Th\'eor\`eme} 1]{Sklar59}, which is commonly referred to as the {\it subcopula} \citep[Definition 3]{Schweizer74}. Defined on $\Ran \FX \times \Ran \FY$, the subcopula must adjust to $\FX$ and $\FY$, and thus {\it cannot} be margin-free -- and neither can any copula subtended by it. In particular, the copula obtained by `jittering' discrete variables with uniform noise to make them artificially continuous, known as  the {\it multilinear} or {\it checkerboard copula} $C^\maltese$, admits a density which amounts to the collection of Pearson ratios such as (\ref{eqn:psi}) \citep[Definition 2.1]{Genest14}, and these were seen not to be variation-independent to the marginal parameters. 
	
	\ppn This may be illustrated further Within the bivariate Bernoulli framework (Example \ref{ex:bivBern}), where (\ref{eqn:Sklar}) reduces down to $p_{00} = C(1-p_{1\bullet},1-p_{\bullet 1})$, for given $p_{1\bullet},p_{\bullet 1} \in (0,1)$ and copula $C$. The value of $C$ at $(1-p_{1\bullet},1-p_{\bullet 1})$ is the subcopula, and is enough to identify the other values $p_{01}$, $p_{10}$ and $p_{11}$ -- hence the whole distribution -- by substitution. The odds ratio, entirely capturing the dependence $\dep$ here, is 
	\[ \omega \doteq \omega(p_{1\bullet},p_{\bullet 1}) = \frac{C(1-p_{1\bullet},1-p_{\bullet 1})(C(1-p_{1\bullet},1-p_{\bullet 1})+p_{1\bullet}+p_{\bullet 1}-1)}{(1-p_{1\bullet} - C(1-p_{1\bullet},1-p_{\bullet 1}))(1-p_{\bullet 1} - C(1-p_{1\bullet},1-p_{\bullet 1}))},\]
	which clearly depends on $p_{1\bullet}$ and $p_{\bullet 1}$ in general -- in fact, the only copula $C$ guaranteeing this function $\omega(p_{1\bullet},p_{\bullet 1})$ to be constant in $p_{1\bullet}$ and $p_{\bullet 1}$ is the {\it Plackett copula}, which was precisely designed in that purpose \citep{Plackett65}. Thus, we may easily construct two bivariate Bernoulli distributions using the {\it same copula} $C$ in (\ref{eqn:Sklar}), but showing very {\it different} dependence structures. \citet[Example 1.1]{Marshall96} shows an extreme example of this. 
	
	\ppn All in all, when both $X$ and $Y$ are continuous variables, then the unique copula $C$ of their distribution is a valid representation of their dependence. In non-continuous cases, on the other hand, neither the subcopula nor any of its copula extension may isolate the dependence structure, so none of these are valid representations of $\dep$.\qed 
\end{example}

\noindent As will be seen, the parities `$\dep$ $\sim$ Pearson's correlation $\rho$', `$\dep$ $\sim$ odds ratio $\omega$', `$\dep$ $\sim$ copula $C$' in the Gaussian, Bernoulli and unspecified continuous models, respectively, are just particular cases of a very general result: there exists a universal representation of the dependence $\dep$ which is valid for any bivariate distribution $\piXY$, and $\rho$, $\omega$ and $C$ will be shown (in Section \ref{subsec:partcases}) to be equivalent to that universal representation in the specific situations of Examples \ref{ex:Gauss}, \ref{ex:bivBern} and \ref{ex:bivCop}.

\section{The dependence $\dep$ and its defining properties} \label{sec:dep}

\subsection{Framework}  \label{subsec:fw}

We consider a general measure-theoretic framework, so as to provide in a unified way results valid across all cases of potential interest (discrete, continuous, mixed or hybrid random vectors alike). Let $(\Xss,\Bs_\Xss,\mu_\Xss)$ and $(\Yss,\Bs_\Yss,\mu_\Yss)$ be two finite measure spaces, and let $\powerset_{\Xss \times \Yss}$ be the set of all probability distributions on $(\Xss \times \Yss,\Bs_\Xss \otimes \Bs_\Yss)$ which are absolutely continuous with respect to $\mu_\Xss \times \mu_\Yss$. For any random vector $(X,Y)$ with distribution $\piXY \in \powerset_{\Xss \times \Yss}$, denote $\phiXY \doteq  \diff \piXY/ \diff (\mu_\Xss \times \mu_\Yss)$. The marginal distributions $\piX$ on $(\Xss,\Bs_\Xss)$ and $\piY$ on $(\Yss,\Bs_\Yss)$ are absolutely continuous with respect to $\mu_\Xss$ and $\mu_\Yss$, respectively, and we let $\phiX \doteq \diff \piX/\diff \mu_\Xss$ and $\phiY \doteq \diff \piY/\diff \mu_\Yss$. We call these functions the joint, $X$-marginal and $Y$-marginal densities of $\piXY$. Also, let the {conditional densities} be 
\begin{equation} \phiYgX=  \frac{\diff \piXY}{\diff (\piX \times \mu_\Yss)}\qquad \text{ and } \qquad \phiXgY  = \frac{\diff \piXY}{\diff (\mu_\Xss \times \piY)}, \label{eqn:conddens} \end{equation}
defined on the sets $\Ss_X \times \Yss$ and $\Xss \times \Ss_Y$, respectively, where $\Ss_X = \{x \in \Xss :\piX(A) > 0 \ \forall A \in \Bs_\Xss \text{ s.t. } A \ni x\} \subseteq \Xss$ and $\Ss_Y = \{y \in \Yss :\piY(B) > 0 \ \forall B \in \Bs_\Yss \text{ s.t. } B \ni y\}  \subseteq \Yss$. Finally, denote $\Ss_{XY} = \{(x,y) \in \Xss  \times \Yss:  \piXY(A \times B) > 0 \ \forall (A,B) \in \Bs_\Xss \otimes \Bs_\Yss \text{ s.t. } A \ni x, B \ni y\} \subseteq \Ss_X \times \Ss_Y \subseteq \Xss \times \Yss$. We will refer to the sets $\Ss_X, \Ss_Y$ and $\Ss_{XY}$ as the (marginal and joint) {supports} of $\piX$, $\piY$ and $\piXY$. The above densities are unique up to relevant almost-everywhere equalities.

\subsection{Conditional densities and dependence} \label{subsec:conddensdep}

It appears from (\ref{eqn:conddens}) that
\begin{align} \piXY(\cdot)  = \iint \phiYgX \indic{\cdot} \diff (\piX \times \mu_\Yss)   = \iint \phiXgY \indic{\cdot}  \diff ( \mu_\Xss \times \piY) \label{eqn:condYX} \end{align}
where $\indic{\cdot}$ is the indicator function of a measurable set. Take two distinct distributions $\pi_{{\scriptscriptstyle XY}}^{(1)}$ and $\pi_{{\scriptscriptstyle XY}}^{(2)}$ in the Fr\'echet class $\Fs(\piX,\piY)$. As $\pi_{{\scriptscriptstyle X}}^{(1)} = \pi_{{\scriptscriptstyle X}}^{(2)} = \piX$, it follows from the first equality that what differentiates $\pi_{{\scriptscriptstyle XY}}^{(1)}$ and $\pi_{{\scriptscriptstyle XY}}^{(2)}$ must be identifiable {only} from contrasting $\{\varphi_{{\scriptscriptstyle Y|X}}^{(1)}(\cdot |x): x\in \Ss_X\}$ and $\{\varphi_{{\scriptscriptstyle Y|X}}^{(2)}(\cdot |x): x\in \Ss_X\}$. It must also be identifiable {only} from contrasting $\{\varphi_{{\scriptscriptstyle X|Y}}^{(1)}(\cdot |y): y\in \Ss_Y\}$ and $\{\varphi_{{\scriptscriptstyle X|Y}}^{(2)}(\cdot |y): y\in \Ss_Y\}$, from the second equality, as $\pi_{{\scriptscriptstyle Y}}^{(1)} = \pi_{{\scriptscriptstyle Y}}^{(2)}  = \piY$. By Definition \ref{def:dep}, what distinguishes $\pi_{{\scriptscriptstyle XY}}^{(1)}$ and $\pi_{{\scriptscriptstyle XY}}^{(2)}$ is their dependence. So:
\begin{fact} \label{fact:conddens} The dependence $\dep$ of a bivariate distribution $\piXY$ must be entirely recoverable from either of its sets of conditional densities $\{\phiYgX\} \doteq \{\phiYgX(\cdot |x): x\in \Ss_X\}$ or $\{\phiXgY\} \doteq \{\phiXgY(\cdot |y):y \in \Ss_Y \}$. We will say that $\dep$ is {\it encapsulated} in $\{\phiYgX\}$ and in $\{\phiXgY\}$.
\end{fact}
\noindent The link between conditional densities and independence/non-independence is well understood. For example, \cite{Wermuth05} describe: `{If proper variables are statistically independent, then the distribution of one of them is the same no matter at which fixed levels the other variable is considered}'; whereas variables not independent are `{such that there are differences in the distributions of one variable for at least some of the levels of the other}'. That is, if the value taken by one of the variables changes, the stochastic behaviour of the other is disturbed -- this is the `influence' mentioned in Section \ref{sec:intro}. In our notation, this amounts to:\footnote{For any measure space $(\Ss,\Bs,\mu)$, we define, for $p>0$, $L_p(\Ss,\mu) \doteq \{\text{measurable }f:\Ss \to \R: \left(\int |f|^p \diff \mu\right)^{1/p} < \infty \}$. Likewise, $L_\infty(\Ss,\mu) \doteq \{\text{measurable }f:\Ss \to \R: \text{ess sup} |f| < \infty$\}.}
\begin{equation*} X \indep Y \iff \exists\, f_{\scriptscriptstyle Y} \in L_1(\Yss,\mu_\Yss) \text{ s.t. } \phiYgX(y|x) = f_{\scriptscriptstyle Y}(y) \qquad (\piX \times \mu_\Yss)\text{-a.e.}, 
\end{equation*} 
i.e., $\phiYgX(y|x)$ is essentially constant in $x$. It is easily shown \citep{Dawid79} that:%\footnote{All statements are meant to hold almost everywhere.}
\begin{align}
	X \indep Y & \iff \exists\, f_{\scriptscriptstyle Y} \in L_1(\Yss,\mu_\Yss) \text{ s.t. } \phiYgX(y|x) = f_{\scriptscriptstyle Y}(y) &  (\piX \times \mu_\Yss)\text{-a.e.} \label{eqn:indep0}  \\ & \iff \phiYgX(y|x) = \phiY(y) &  (\piX \times \mu_\Yss)\text{-a.e.}  \notag \\
	& \iff \phiXY(x,y) = \phiX(x) \phiY(y)  &  (\mu_\Xss \times \mu_\Yss)\text{-a.e.}   \label{eqn:indep2}\\
	& \iff \phiXgY(x|y) = \phiX(x) &  (\mu_\Xss \times \piY)\text{-a.e.}  \notag  \\ & \iff  \exists\, f_{\scriptscriptstyle X} \in L_1(\Xss,\mu_\Xss) \text{ s.t. }  \phiXgY(x|y) =  f_{\scriptscriptstyle X}(x) &  (\mu_\Xss \times \piY)\text{-a.e.}. \label{eqn:indep3}
\end{align}
From (\ref{eqn:indep0}) and (\ref{eqn:indep3}), it is clear that (non-)independence is determined exclusively by the conditional densities, without any reference to the marginals. Fact \ref{fact:conddens} actually strengthens this statement, showing that {\it full characterisation} of the dependence $\dep$ of $\piXY$ is possible from either one of $\{\phiYgX\}$ or $\{\phiXgY\}$, and only from it. It follows that:

\begin{corollary} \label{cor:conddens} Two distributions $\pi_{{\scriptscriptstyle XY}}^{(1)}$ and $\pi_{{\scriptscriptstyle XY}}^{(2)}$ sharing the same $X|Y$-conditional density (i.e, $\{\varphi_{{\scriptscriptstyle X|Y}}^{(1)}\} \equiv \{\varphi_{{\scriptscriptstyle X|Y}}^{(2)}\}$) \uline{or} the same $Y|X$-conditional density (i.e, $\{\varphi_{{\scriptscriptstyle Y|X}}^{(1)}\} \equiv \{\varphi_{{\scriptscriptstyle Y|X}}^{(2)}\}$) have identical dependence structures ($\dep^{(1)} = \dep^{(2)}$).
\end{corollary}

\subsection{One-to-one marginal transformations and canonical supports} \label{subsec:equiv}

An immediate consequence of Fact \ref{fact:conddens} is that $\dep$ is essentially invariant under one-to-one `{marginal transformations}'. Formally, let $(\Xss',\Bs_{\Xss'})$ be a measurable space isomorphic to $(\Xss,\Bs_{\Xss})$ \cite[Chapter 13]{Dudley02}. For a given one-to-one function $\Phi:\Xss \to \Xss'$, such that both $\Phi$ and $\Phi^{-1}$ are measurable (i.e., $\Phi$ is {bimeasurable}), define the random variable $X' \doteq \Phi(X)$, with distribution $\pi_{\scriptscriptstyle X'} = \piX \circ \Phi^{-1}$ (`{push-forward}') on the support $\Ss_{X'} = \{x' \in \Xss' :\pi_{\scriptscriptstyle X'}(A') > 0 \ \forall A' \in \Bs_{\Xss'} \text{ s.t. } A' \ni x'\}= \{x' \in \Xss' :\piX(A) > 0 \ \forall A \in \Bs_{\Xss} \text{ s.t. } A  \ni \Phi^{-1}(x')\}$. The distribution $\pi_{\scriptscriptstyle X' Y}$ of the vector $(X',Y)$ admits as conditional density $\{\varphi_{\scriptscriptstyle Y|X'}\} = \{ \varphi_{\scriptscriptstyle Y|X'}(\cdot|x'):x'\in \Ss_{X'}\}$, with
\[\varphi_{\scriptscriptstyle Y|X'}(\cdot |x')  = \varphi_{\scriptscriptstyle Y|\Phi^{-1}(X')}(\cdot|\Phi^{-1}(x')) = \varphi_{\scriptscriptstyle Y|X}(\cdot |\Phi^{-1}(x')).\]
From a conditioning-on-$X$ perspective, the transformation $X' = \Phi(X)$ only changes the `label' assigned to the value on which we condition, but leaves the effective conditional density of $Y$ unaffected. Similarly to (\ref{eqn:condYX}) we can write
\[\piXY(\cdot)  = \iint \varphi_{\scriptscriptstyle Y|X'}(y|\Phi(x)) \indic{\cdot} \diff (\piX \times\mu_\Yss) (x,y),  \]
showing that $\dep$ would be entirely recoverable from $\varphi_{\scriptscriptstyle Y|X'}$, too, provided we know $\Phi$. Also, as $\dep$ provides all required information for unequivocally identifying $\piXY$ in $\Fs(\piX,\piY)$, it inherently allows unequivocal identification of $\pi_{\scriptscriptstyle X' Y}$ in $\Fs(\pi_{\scriptscriptstyle X'},\piY)$ as well, as distributions of these two Fr\'echet classes are in one-to-one correspondence through the given $\Phi$. 

\ppn We can replicate this with a bimeasurable bijection $\Psi: \Yss \to \Yss'$ which defines a random variable $Y' = \Psi(Y)$ living on another measurable space $(\Yss',\Bs_{\Yss'})$. Consequently, for any given mappings $\Phi$ and $\Psi$ as above, the respective dependence structures $\dep$ and $\dep'$ of the vectors $(X,Y)$ and $(X',Y')$ are equivalent (one can be recovered from the other). Thus, any valid representation of $\dep$ must be {equivariant} \cite[Chapter 3]{Lehmann06} under such one-to-one marginal transformations: the information they contain is effectively the same, it just need be encoded differently so as to keep track of the `changes of label' $x'=\Phi(x)$ and $y'=\Psi(y)$. 

\ppn It follows that the dependence structures of all bivariate distributions with {\it marginal supports of the same cardinality} (the mappings $\Phi$ and $\Psi$ are one-to-one) are objects of the same kind. Therefore, dependence structures of distributions with marginals supported on arbitrary `canonical' sets $\Ss'_X \subseteq \Xss'$ and $\Ss'_Y \subseteq \Yss'$ are enough for characterising the dependence between marginals supported on any sets $\Ss_{X} \subseteq \Xss$ and $\Ss_{Y} \subseteq \Yss$ such that $|\Ss_{X}|=|\Ss'_X|$ and $|\Ss_{Y}|=|\Ss'_Y|$.\footnote{$|\Ss|$ denotes the cardinality of the set $\Ss$.} E.g., the dependence between discrete variables supported on sets $\Ss_{X}$ and $\Ss_{Y}$ of finite cardinality $R$ and $S$ can be entirely described by the dependence of a distribution with margins supported on $\Ss'_X = \{0,\ldots,R-1\}$ and $\Ss'_Y = \{0,\ldots,S-1\}$. Likewise, dependence structures between continuous variables both supported on $\Ss'_X = \Ss'_Y = [0,1]$ are sufficient for describing dependence in any bivariate continuous distribution. We shall call versions of $\piXY$ with marginal supports adjusted to the agreed canonical sets of the right cardinality, `{\it canonical versions}' of $\piXY$.

\subsection{Joint support and regional dependence} \label{subsec:regio}

Clearly, $X$ and $Y$ can only be independent if $\Ss_{XY} = \Ss_{X} \times \Ss_{Y}$ (`{rectangular support}'): if not, there exist sets $A \in \Bs_\Xss$ and $B \in\Bs_\Yss$ of positive $\piX$- and $\piY$-measure such that $\phiXY(x,y)=0$ for almost all $(x,y) \in A \times B$, violating (\ref{eqn:indep2}). Thus, any non-rectangular support $\Ss_{XY}$ automatically implies non-independence, through the fact that some values of $X$ and $Y$ are not allowed simultaneously -- a strong form of mutual `influence', indeed. \cite{Holland86} refer to this as `{\it regional dependence}'. 

\ppn Regional dependence is dictated by the {shape} of $\Ss_{XY}$, but not by the actual distribution $\piXY$ supported on that set. As the shape of $\Ss_{XY}$ cannot be recovered from marginal specifications only, regional dependence is part of $\dep$, as per Definition \ref{def:dep}. The {shape} of $\Ss_{XY}$ should be distinguished from the set $\Ss_{XY}$ itself, though. In particular, $\Ss_{XY}$ contains marginal information, such as the supports $\Ss_X$ and $\Ss_Y$, so cannot be entirely part of $\dep$. This said, the shape of $\Ss_{XY}$ can be inferred from the joint support $\Ss'_{XY}$ of any canonical version $\pi'_{\scriptscriptstyle XY}$ of $\piXY$ (Section \ref{subsec:equiv}) -- and this one can only reveal the uninformative $\Ss_X'$ and $\Ss_Y'$. Thus, to comply with Definition \ref{def:dep}, the dependence $\dep$ is always to be described in terms of canonical versions of distributions. 

\ppn Therefore, all distributions coming into play below will be assumed to be `canonical', meaning that their marginal supports have been preliminarily transformed to some universally agreed sets of the right cardinality.\footnote{Their exact specification is irrelevant, as they are arbitrary.} Thus, it will always be assumed that the joint support of a distribution does not contain any crucial information about its marginals (as $\Ss_X,\Ss_Y$ are universally known), and that marginal supports of the same cardinality are equal; that is, for any two distributions $\pi^{(1)}_{\scriptscriptstyle XY}$ and $\pi^{(2)}_{\scriptscriptstyle XY}$, $|\Ss^{(1)}_X| = |\Ss^{(2)}_X| \Rightarrow  \Ss^{(1)}_X = \Ss^{(2)}_X$ and $|\Ss^{(1)}_Y| = |\Ss^{(2)}_Y| \Rightarrow  \Ss^{(1)}_Y = \Ss^{(2)}_Y$. In particular, as any two distributions $\pi^{(1)}_{\scriptscriptstyle XY},\pi^{(2)}_{\scriptscriptstyle XY}$ sharing the same dependence automatically share the same regional dependence, this assumption allows us to write:
\begin{equation} \dep^{(1)} = \dep^{(2)} \Rightarrow \Ss^{(1)}_{XY} = \Ss^{(2)}_{XY}.\label{eqn:depsup} \end{equation}
All this understood, no further explicit reference to the canonical versions/supports will be made.

\subsection{Marginal replacements} \label{subsec:margrep}

Another consequence of Fact \ref{fact:conddens} is that $\dep$ must be invariant to `{marginal replacement}' \citep[Definition 4.1]{Holland87}. Let $\piXY \in \Fs(\piX,\piY)$ and consider another $X$-marginal measure $\piXstar$ equivalent to $\piX$ ($\piX \ll \piXstar \ll \piX$, denoted here $\piX\approx \piXstar$). One can use the conditional density $\{\phiYgX\}$ of $\piXY$ as a Markov kernel, and construct the distribution $\piXYstar$ with density
\begin{equation} \phiXYstar = \phiYgX \frac{\diff \piXstar}{\diff \mu_\Xss} = \phiXY \frac{\diff \piXstar}{\diff \piX} . \label{eqn:pistar} \end{equation}
Although $\piXstar \neq  \piX$ and $\piYstar(\cdot) = \iint  \phiXY \frac{\diff \piXstar}{\diff \piX} \indic{\Ss_X \times \cdot}\diff (\mu_\Xss \times \mu_\Yss) \neq  \iint  \phiXY \indic{\Ss_X \times \cdot} \diff (\mu_\Xss \times \mu_\Yss) = \piY(\cdot)$ in general, the dependence structure of $\piXY$ and $\piXYstar$ is identical as they share the same $\{\phiYgX\}$ (Corollary \ref{cor:conddens}). Re-writing (\ref{eqn:pistar}) as 
\begin{equation} \frac{\diff \piXYstar}{\diff \piXY} = \frac{\diff \piXstar}{\diff \piX} \label{eqn:pistar2} \end{equation}
shows that $\piXYstar \approx \piXY$ (as $\piXstar \approx \piX$) and thus $\piYstar \approx \piY$. From any other $Y$-marginal measure $\pi^{**}_{\scriptscriptstyle Y} \approx \piY$, one can produce by the same token yet another joint distribution $\pi^{**}_{\scriptscriptstyle XY}$ with density
\begin{equation} \varphi^{**}_{\scriptscriptstyle XY} = \phiXgYstar \frac{\diff \pi^{**}_{\scriptscriptstyle Y}}{\diff \mu_\Yss} = \phiXYstar \frac{\diff \pi^{**}_{\scriptscriptstyle Y}}{\diff \piYstar} = \phiXY \frac{\diff \piXstar}{\diff \piX}\frac{\diff \pi^{**}_{\scriptscriptstyle Y}}{\diff \piYstar}. \label{eqn:pistarstar} \end{equation}
Hence
\begin{equation} \frac{\diff \pi_{\scriptscriptstyle XY}^{**}}{\diff \piXY} = \frac{\diff \piXstar}{\diff \piX} \frac{\diff \pi^{**}_{\scriptscriptstyle Y}}{\diff \piYstar}, \label{eqn:IPF1step} \end{equation}
showing that $\pi_{\scriptscriptstyle XY}^{**} \approx \piXY$. This time  $\{\varphi^{**}_{\scriptscriptstyle X|Y}\} = \{\phiXgYstar\}$, and $\pi_{\scriptscriptstyle XY}^{**}$ has the same dependence as $\piXYstar$ (Corollary \ref{cor:conddens}), and hence as $\piXY$ too. We have thus created two distributions $\piXYstar$ and $\pi^{**}_{\scriptscriptstyle XY}$ with respective $X$- and $Y$-marginals set to arbitrary $\piXstar \approx \piX$ and $\pi^{**}_{\scriptscriptstyle Y} \approx \piY$ while maintaining the same dependence $\dep$ as the initial $\piXY$. Note that $\pi^{**}_{\scriptscriptstyle XY}$ does not share any of the conditional densities $\{\phiYgX\}$ or $\{\phiXgY\}$ of $\piXY$, illustrating that Corollary \ref{cor:conddens} provides a sufficient but not necessary condition for equal dependence structures. In (\ref{eqn:pistar2}) and (\ref{eqn:IPF1step}), the Radon-Nikodym derivatives of $\pi^{*}_{\scriptscriptstyle XY}$ and $\pi^{**}_{\scriptscriptstyle XY}$ with respect to $\piXY$ have a product form -- Theorem \ref{thm:equivclass} below establishes that it is not coincidental. Finally, as $\pi_{\scriptscriptstyle XY}^{**} \approx \pi_{\scriptscriptstyle XY}^{*} \approx \piXY$, the joint support $\Ss_{XY}$ is preserved through the successive replacements, in agreement with (\ref{eqn:depsup}).

\subsection{Iterative Proportional Fitting procedure and $\Is$-projections} \label{sec:IPF}

The marginal replacement of $\piX$ by $\piXstar$ via (\ref{eqn:pistar}) affects the $Y$-margin (it was noted that $\piYstar \neq \piY$). Likewise, the marginal replacement of $\piYstar$ by $\pi_{\scriptscriptstyle Y}^{**}$ in (\ref{eqn:pistarstar}) affects the $X$-margin (in particular: $\pi^{**}_{\scriptscriptstyle X} \neq \piXstar$). Thus, for a couple of targeted marginal distributions $(\piXstar,\piYstar)$, it is not clear if we can obtain from $\piXY$ a distribution of $\Fs(\piXstar,\piYstar)$ via a similar dependence-preserving process. The following lemma gives a first element of answer. Let $\piXY \in \powerset_{\Xss \times \Yss}$, and define 
\begin{equation*} \text{N}_{\piXY} = \{(A,B) \in \Bs_\Xss \otimes \Bs_\Yss \text{ s.t.\ } \piXY(A \times B) = 0 \}, \label{eqn:NXY} \end{equation*}
the collection of rectangular subsets of $\Xss \times \Yss$ on which $\piXY$ is null. Also, let $\Sigma_{\piXY} = \{\pi \in \powerset_{\Xss \times \Yss}: \pi \ll \piXY \}$, that is the distributions of $\powerset_{\Xss \times \Yss}$ with support contained entirely in $\Ss_{XY}$. 

\begin{lemma} \label{lem:Leuridan} Let $\piXY \in \powerset_{\Xss \times \Yss} \cap \Fs(\piX,\piY)$ have support $\Ss_{XY}$. For any couple of marginal distributions $(\piXstar,\piYstar)$ such that $\piXstar \approx \piX$ and $\piYstar \approx \piY$, we have:
	\begin{enumerate}
		\item \label{item:NSC} $\powerset_{\Xss \times \Yss} \cap \Sigma_{\piXY} \cap \Fs(\piXstar,\piYstar) \neq \emptyset$ $\iff$  $\forall\,(A,B) \in \textup{N}_{\piXY}, \piXstar(A) + \piYstar(B) \leq  1$; 
		\item \label{item:NSC2} If $\powerset_{\Xss \times \Yss} \cap \Sigma_{\piXY} \cap \Fs(\piXstar,\piYstar) \neq \emptyset$:
		\begin{enumerate}
			\item \label{item:exist} $\exists\  \overline{\pi}^*_{\scriptscriptstyle XY} \in \powerset_{\Xss \times \Yss} \cap \Sigma_{\piXY} \cap \Fs(\piXstar,\piYstar)$ such that $\piXYstar \in \powerset_{\Xss \times \Yss} \cap \Sigma_{\piXY} \cap \Fs(\piXstar,\piYstar)$ $\Rightarrow$ $\piXYstar \ll \overline{\pi}^*_{\scriptscriptstyle XY}$;
			\item \label{item:exist2} $\overline{\pi}^*_{\scriptscriptstyle XY} \approx \piXY$ $\iff$ $\left[ \forall(A,B) \in \textup{N}_{\piXY},\  \piXstar(A) + \piYstar(B) = 1 \Rightarrow (A^c,B^c) \in \textup{N}_{\piXY} \right]$, where $A^c$ and $B^c$ are the complements of $A$ and $B$ in $\Xss$ and $\Yss$, respectively.
		\end{enumerate}
	\end{enumerate}
\end{lemma}

\noindent  Part \ref{item:NSC} gives a necessary and sufficient condition for the existence of some distribution of $\Fs(\piXstar,\piYstar)$ with support contained in $\Ss_{XY}$. If so, Part \ref{item:NSC2}-\ref{item:exist}) says that there exists one of them whose support contains the support of all the others. Part \ref{item:NSC2}-\ref{item:exist2}) gives a necessary and sufficient condition for that `maximal support' to be $\Ss_{XY}$ -- in other words, for the existence of a distribution of $\Fs(\piXstar,\piYstar)$ with support $\Ss_{XY}$. We conclude:

\begin{corollary} \label{cor:Leuridan2} Let $\piXY \in \powerset_{\Xss \times \Yss} \cap \Fs(\piX,\piY)$. For any couple of marginal distributions $(\piXstar,\piYstar)$ such that $\piXstar \approx \piX$ and $\piYstar \approx \piY$, there exists a distribution of $\piXYstar \in \Fs(\piXstar,\piYstar)$ such that $\piXYstar \approx \piXY$ if and only if 
	\begin{equation}  \begin{array}{l} \forall (A,B) \in \textup{N}_{\piXY},\ 	 \piXstar(A) + \piYstar(B) \leq 1, \text{ and } \\  \qquad \text{ any } (A,B) \in \textup{N}_{\piXY} \text{ s.t. } \piXstar(A) + \piYstar(B) = 1  \text{ is such that }(A^c,B^c) \in \textup{N}_{\piXY}.  \end{array}\label{eqn:compacond} \end{equation}
\end{corollary}

\noindent The importance of Corollary \ref{cor:Leuridan2} stems from (\ref{eqn:depsup}), establishing the connection between the joint support of a distribution and its dependence. Hence (\ref{eqn:compacond}) provides a {necessary} condition for finding a distribution $\piXYstar$ in $\Fs(\piXstar,\piYstar)$ with the same dependence as $\piXY$: marginals $(\piXstar,\piYstar)$ violating (\ref{eqn:compacond}) cannot co-exist with the support of $\piXY$, that is, with its specific regional dependence. Condition (\ref{eqn:compacond}) is trivially fulfilled for all $(\piXstar,\piYstar)$ if $\piXY$ has rectangular support, as in that case any $(A,B) \in \textup{N}_{\piXY}$ is such that $\piX(A) = \piY(B) = 0$, thus $\piXstar(A) = \piYstar(B) = 0$ as $\piXstar \ll \piX$ and $\piYstar \ll \piY$. In fact, (\ref{eqn:compacond}) is a rather mild condition, as breaking it requires the existence of a `large' area in $\Ss_X \times \Ss_Y$ on which $\piXY$ is null, imposing extensive incompatibility between values of $X$ and $Y$ (`{extreme}' regional dependence). 

\ppn An example may be the case of a bivariate Bernoulli distribution (\ref{eqn:tab22}) with the main diagonal as support, that is, $p_{10} = p_{01} = 0$. This evidently implies $p_{1\bullet} = p_{\bullet 1}$, that is, identical margins $\piX= \piY$. Effectively this specific regional dependence structure imposes $X=Y$, therefore it cannot be found in Fr\'echet classes $\Fs(\text{Bern}(p_{1\bullet}),\text{Bern}(p_{\bullet 1 }))$ such that $p_{1\bullet} \neq p_{\bullet 1}$, for the obvious reason. 

\ppn We show below that (\ref{eqn:compacond}) is actually also {sufficient} for finding a distribution $\piXYstar$ in $\Fs(\piXstar,\piYstar)$ with the same dependence as $\piXY$, with a few added mild technical conditions on $\piXY$.\footnote{These conditions are those given in \cite{Nussbaum93}. In fact, on page 81, \cite{Nussbaum93} mentioned that these extra assumptions may well be an artefact of the argument, and the result of main interest [his Theorem 3.9] may remain true without them -- this author has not attempted to prove or disprove that claim, arguing in Appendix \ref{app:A} that these requirements are, in any case, all but constraining in the considered framework.}  The subset of distributions of $\powerset_{\Xss \times \Yss}$ satisfying those conditions, listed in Appendix \ref{app:A}, will be denoted $\overline{\powerset}_{\Xss \times \Yss}$. First, we need the following lemma, which easily ensues from \cite{Nussbaum93}'s Theorem 3.9, Lemma 3.11 and Corollary 3.13 (his `compatibility' condition amounts to (\ref{eqn:compacond}), while his Hypotheses 3.6-3.10 correspond to $\piXY \in \overline{\powerset}_{\Xss \times \Yss}$).

\begin{lemma} \label{lem:nussbaum} Let $\piXY \in \overline{\powerset}_{\Xss \times \Yss} \cap \Fs(\piX,\piY)$ and $(\piXstar,\piYstar)$ be a couple of marginal distributions such that $\piXstar \approx \piX$ and $\piYstar \approx \piY$. Under condition (\ref{eqn:compacond}), there exists a unique distribution $\piXYstar \in \Fs(\piXstar,\piYstar)$ such that
	\begin{equation} \frac{\diff \piXYstar}{\diff \piXY} = {\alpha}^*_{\scriptscriptstyle X} {\beta}^*_{\scriptscriptstyle Y} \label{eqn:alphabeta}\end{equation}
	for (a.e.-)positive functions ${\alpha}^*_{\scriptscriptstyle X} \in L_1(\Xss,\mu_\Xss)$ and ${\beta}^*_{\scriptscriptstyle Y} \in L_1(\Yss,\mu_\Yss)$ unique up to a multiplicative constant.
\end{lemma}

\noindent According to \citet[Theorem 4.19 and Remark 4.20]{Nussbaum93}, the distribution $\piXYstar$ in (\ref{eqn:alphabeta}) can be obtained as the limit of a (geometrically convergent) iterative procedure  consisting in successive $X$- and $Y$-marginal replacements {\it \`a la} (\ref{eqn:pistar})-(\ref{eqn:pistarstar}) alternately. Specifically, if one sets
\begin{equation} \arraycolsep=1.4pt\def\arraystretch{1.5} \forall n \geq 0, \quad \left\{\begin{array}{l l} \frac{\diff \pi^{(2n)}_{\scriptscriptstyle XY}}{\diff \piXY} = \alpha^{(n)}_{\scriptscriptstyle X} \beta^{(n)}_{\scriptscriptstyle Y}, & \frac{\diff \pi^{(2n+1)}_{\scriptscriptstyle XY}}{\diff \piXY} = \alpha^{(n+1)}_{\scriptscriptstyle X} \beta^{(n)}_{\scriptscriptstyle Y} \\ \alpha^{(n+1)}_{\scriptscriptstyle X} = \frac{\diff \piXstar}{\diff \pi^{(2n)}_{\scriptscriptstyle X}}\alpha^{(n)}_{\scriptscriptstyle X} , & \beta^{(n+1)}_{\scriptscriptstyle Y} = \frac{\diff \piYstar}{\diff \pi^{(2n+1)}_{\scriptscriptstyle Y}}\beta^{(n)}_{\scriptscriptstyle Y} \end{array}\right., \label{eqn:IPF1}\end{equation}
with $\alpha^{(0)}_{\scriptscriptstyle X} = \beta^{(0)}_{\scriptscriptstyle Y} =1$ and $\pi^{(0)}_{\scriptscriptstyle X} = \piX$, $\pi^{(0)}_{\scriptscriptstyle Y} = \piY$ (and $\forall n' \geq 0$, $\pi^{(n')}_{\scriptscriptstyle X},\pi^{(n')}_{\scriptscriptstyle Y}$ are the margins of $\pi^{(n')}_{\scriptscriptstyle XY}$), then $\pi^{(n')}_{\scriptscriptstyle XY} \to \piXYstar$ in total variation as $n'\to\infty$. Given that each step corresponds to a dependence-preserving marginal replacement and that $\pi^{(n')}_{\scriptscriptstyle XY} \approx \piXY\ \forall n'$ including in the limit ($\piXYstar \approx \piXY$), the initial dependence of $\piXY$ is maintained throughout the process. Thus:

\begin{corollary} \label{cor:alphabetadep} Let $\piXY \in \overline{\powerset}_{\Xss \times \Yss} \cap \Fs(\piX,\piY)$ and $(\piXstar,\piYstar)$ be marginal distributions such that $\piXstar \approx \piX$ and $\piYstar \approx \piY$. Under condition (\ref{eqn:compacond}), there exists a unique distribution $\piXYstar$ in $\Fs(\piXstar,\piYstar)$ sharing the same dependence $\dep$ as $\piXY$; and that distribution satisfies (\ref{eqn:alphabeta}).
\end{corollary}

\noindent The recursive process (\ref{eqn:IPF1}) is known as the {Iterative Proportional Fitting} (IPF) procedure \citep{Deming40,Sinkhorn64,Ireland68,Kullback68,Fienberg70,Zalovznik11,Idel16}. It has strong links to the so-called `$\Is$-geometry' \citep{Csiszar75}. For two distributions $\widetilde{\pi},\pi$ on $(\Xss \times \Yss,\Bs_\Xss \otimes \Bs_\Yss)$, define 
\begin{equation*} \Is(\widetilde{\pi}\| \pi) = \left\{\begin{array}{l l} \iint \log \left(\frac{\diff \widetilde{\pi}}{\diff \pi} \right)\diff \widetilde{\pi} = \iint \log \left(\frac{\diff \widetilde{\pi}}{\diff \pi} \right)\frac{\diff \widetilde{\pi}}{\diff \pi}  \diff \pi & \text{ if } \widetilde{\pi} \ll \pi \\ \infty & \text{ if }  \widetilde{\pi} \not\ll \pi\end{array}\right., \label{eqn:KL2}\end{equation*}
the relative entropy of $\widetilde{\pi}$ with respect to $\pi$ -- also called the $\Is$-divergence, or Kullback-Leibler divergence. Properties of this quantity are well-known \citep{Cover06}. In particular, $\Is(\widetilde{\pi}\| \pi) \geq 0$ always, with $\Is(\widetilde{\pi}\| \pi) = 0 \iff \pi = \widetilde{\pi}$. Hence, although it is not symmetric (and does not satisfy a triangular inequality), $\Is(\cdot\|\cdot)$ may be thought of as some sort of distance between distributions, which defines a specific geometry. {\it Inter alia}, we may define a concept of projection as follows:

\begin{definition} Let $\Es$ be some set of distributions on $(\Xss \times \Yss,\Bs_\Xss \otimes \Bs_\Yss)$, assumed to be convex and closed in total variation. Then, for a distribution $\pi$, any distribution $\widetilde{\pi} \in \Es$ such that $\Is(\widetilde{\pi} \| \pi) = \inf_{\check{\pi} \in \Es} \Is(\check{\pi}	\|\pi)$ is called an $\Is$-{projection} of $\pi$ on $\Es$.
\end{definition}
\noindent An $\Is$-projection exists and is unique as soon as $\inf_{\check{\pi} \in \Es} \Is(\check{\pi}	\|\pi) <\infty$, an obvious `reachability' condition \citep{Csiszar75}. Noting that a Fr\'echet class of distributions is convex and closed in total variation \citep[Lemma 4.1]{Cramer00}, we can define the $\Is$-projection of $\piXY$ on $\Fs(\piXstar,\piYstar)$ as the distribution $\widetilde{\pi}^*_{\scriptscriptstyle XY}$ satisfying
\[\Is(\widetilde{\pi}^*_{\scriptscriptstyle XY}\| \piXY) = \inf_{\pi \in \Fs(\piXstar,\piYstar)} \Is(\pi\|\piXY).  \] 
When $\inf_{\pi \in \Fs(\piXstar,\piYstar)} \Is(\pi\|\piXY)<\infty$ (i.e., when $\widetilde{\pi}^*_{\scriptscriptstyle XY}$ exists and is unique), \cite{Ruschendorf95} showed that the IPF procedure (\ref{eqn:IPF1}) converges to $\widetilde{\pi}^*_{\scriptscriptstyle XY}$. From \citet[Corollary 3.2]{Csiszar75} and \citet[Theorem 3]{Ruschendorf93}, we can further deduce that, under condition (\ref{eqn:compacond}), $\widetilde{\pi}^*_{\scriptscriptstyle XY}$ is in fact $\piXYstar$, the distribution defined by (\ref{eqn:alphabeta}); that is, the distribution of $\Fs(\piXstar,\piYstar)$ with the same dependence as $\piXY$ (Corollary \ref{cor:alphabetadep}). Hence, {\it $\Is$-projections on Fr\'echet classes maintain dependence}.

\ppn Therefore, the $\Is$-geometry seems particularly well suited for providing a formal geometrical framework around (\ref{eqn:decomp}). For example, \cite{Csiszar75} proved a Pythagoras-like theorem for $\Is$: for two arbitrary distributions $\piXY \in \Fs(\piX,\piY)$ and $\check{\pi}^*_{\scriptscriptstyle XY} \in \Fs(\piXstar,\piYstar)$, we have
\[\Is(\check{\pi}^*_{\scriptscriptstyle XY} \| \piXY) = \Is(\check{\pi}^*_{\scriptscriptstyle XY} \| \widetilde{\pi}^*_{\scriptscriptstyle XY}) + \Is(\widetilde{\pi}^*_{\scriptscriptstyle XY} \| \piXY), \]
where $\widetilde{\pi}^*_{\scriptscriptstyle XY}$ is the $\Is$-projection of $\piXY$ on $\Fs(\piXstar,\piYstar)$. Hence $\widetilde{\pi}^*_{\scriptscriptstyle XY}$ has the same marginals as $\check{\pi}^*_{\scriptscriptstyle XY}$ (but not the same dependence) while, from the previous argument, $\widetilde{\pi}^*_{\scriptscriptstyle XY}$ has the same dependence as $\piXY$ (but not the same marginals). Thus the `distance' from $\piXY$ to $\check{\pi}^*_{\scriptscriptstyle XY}$ can be broken down into two {orthogonal} contributions: what it takes to adjust the marginals (maintaining the dependence) plus what it takes to adjust the dependence (maintaining the marginals).

\subsection{Compatibility of conditional densities} \label{sec:compatibility}

Fact \ref{fact:conddens} establishes that full description of $\dep$ may be obtained from either $\{\phiYgX\}$ or $\{\phiXgY\}$. Therefore we may posit the existence of a mapping, say $\Delta$, taking as argument conditional densities such as $\{\phiYgX\}$ or $\{\phiXgY\}$, and returning (a complete representation of) the dependence $\dep$ which they encapsulate. Informally, we would have, for any specification of a $X|Y$-conditional density $\{\phiXgY\}$:
\[\Delta: \{\phiXgY\} \to \Delta(\{\phiXgY\}) = \text{dependence } \dep \text{ encapsulated in }\{\phiXgY\},\]
and for any specification of a $Y|X$-conditional density $\{\phiYgX\}$:
\[\Delta: \{\phiYgX\} \to \Delta(\{\phiYgX\}) = \text{dependence } \dep \text{ encapsulated in }\{\phiYgX\}.\]
That map should necessarily be such that
\begin{equation} \Delta(\{\phiXgY\}) = \Delta(\{\phiYgX\})  \label{eqn:deltas} \end{equation}	
when $\{\phiXgY\}$ and $\{\phiYgX\}$ are the two `twin' sets of conditional densities of a common bivariate distribution $\piXY$, as they encapsulate the same dependence -- that of $\piXY$. 

\ppn From a modelling perspective, it is sometimes more natural to specify an appropriate model for a vector $(X,Y)$ through its two sets of conditional distributions ($X$ given $Y$, and $Y$ given $X$), rather than through a joint distribution directly \citep{Arnold01}. Yet, two separately postulated conditional models, say $\{\phiXgYstar\}$ and $\{\phiYgX\}$, do not necessarily interlock into a proper bivariate distribution. When there exists a bivariate distribution $\widetilde{\pi}_{\scriptscriptstyle XY}$ with $\{\phiXgYstar\}$ and $\{\phiYgX\}$ as conditional densities, these two are called {compatible} \citep{Arnold89,Arnold22,Wang08}. Here we establish that two sets of conditional densities are compatible if and only if they encapsulate the same dependence.

\begin{proposition} \label{prop:compa} For $\Ss_X \subset \Xss$ and $\Ss^*_Y \subset \Yss$, let $\{\phiYgX\} = \{\phiYgX(\cdot |x): x\in \Ss_X\}$ and $\{\phiXgYstar\} = \{\phiXgYstar(\cdot |y):y \in \Ss^*_Y \}$ be some $Y|X$- and a $X|Y$-conditional density specifications. Then $\{\phiYgX\}$ and $\{\phiXgYstar\}$ are two compatible sets of conditional densities if and only if $\{\phiYgX\}$ and $\{\phiXgYstar\}$ both encapsulate the same dependence structure; i.e., $\Delta(\{\phiYgX\}) = \Delta(\{\phiXgYstar\})$.
\end{proposition}

\noindent Intuitively, any arbitrary $X$-marginal distribution $\widetilde{\pi}_{\scriptscriptstyle X}$ and $Y|X$-conditional density $\phiYgX$ can be freely combined as in (\ref{eqn:pistar}) to form a bivariate distribution -- provided that the support of $\widetilde{\pi}_{\scriptscriptstyle X}$ matches $\Ss_X$ in $\{\phiYgX\}$. Likewise, any arbitrary $Y$-marginal and $X|Y$-conditional density can be freely combined, up to the same support restriction. Thus, if we set the putative $\widetilde{\pi}_{\scriptscriptstyle XY} = (\widetilde{\pi}_{\scriptscriptstyle X},\widetilde{\pi}_{\scriptscriptstyle Y};\widetilde{\dep})$ as in (\ref{eqn:decomp}), $\{\phiYgX\}$ does not restrict $\widetilde{\pi}_{\scriptscriptstyle X}$ while $\{\phiXgYstar\}$ does not restrict $\widetilde{\pi}_{\scriptscriptstyle Y}$, and (in)compatibility of $\{\phiYgX\}$ and $\{\phiXgYstar\}$ may only be dictated by their encapsulated dependence $\dep$ and $\dep^*$. Clearly $\dep = \dep^*$  (which implies matching supports, as per (\ref{eqn:depsup})) is all what is needed for completing the representation $\widetilde{\pi}_{\scriptscriptstyle XY} = (\widetilde{\pi}_{\scriptscriptstyle X},\widetilde{\pi}_{\scriptscriptstyle Y};\widetilde{\dep})$ with $\widetilde{\dep} =\dep = \dep^*$. On the other hand, if $\dep \neq \dep^*$, then $\{\phiYgX\}$ and $\{\phiXgYstar\}$ have no common ground around which articulate the desired bivariate distribution. 

\section{Universal representation of the dependence $\dep$} \label{sec:universal0}

\subsection{Equivalence classes of dependence structure} \label{subsec:equivclass} 

Combining the results of Section \ref{sec:dep}, we can now establish:

\begin{theorem} \label{thm:equivclass} Let $\piXY$ and $\piXYstar$ be two distributions in $\overline{\powerset}_{\Xss \times \Yss}$. Then $\piXY$ and $\piXYstar$ admit the same dependence $\dep = \dep^*$ if and only if there exist (a.e.-)positive functions $\alpha^*_{\scriptscriptstyle X} \in L_1(\Xss,\piX)$ and $\beta^*_{\scriptscriptstyle Y} \in L_1(\Yss,\piY)$ such that 
	\begin{equation} \frac{\diff \piXYstar}{\diff \piXY} = \alpha^*_{\scriptscriptstyle X} \beta^*_{\scriptscriptstyle Y}. \label{eqn:equiv} \end{equation}
\end{theorem}

\ppn The Radon-Nikodym derivative (\ref{eqn:equiv}) is (a.e.-)positive, hence $\piXY$ and $\piXYstar$ have the same support -- in agreement with (\ref{eqn:depsup}). The product form  $\alpha^*_{\scriptscriptstyle X} \beta^*_{\scriptscriptstyle Y}$, reminiscent of (\ref{eqn:indep2}), implies that the relative difference between $\piXY$ and $\piXYstar$ shows up as two independent marginal adjustments, unable to affect the inner dependence structure. Evidently (\ref{eqn:equiv}) defines an equivalence relationship between distributions, which in turn defines equivalence classes of dependence:
\begin{equation}   \piXY \sim \piXYstar  \iff \dep = \dep^* \iff  \frac{\diff \piXYstar}{\diff \piXY} = \alpha^*_{\scriptscriptstyle X} \beta^*_{\scriptscriptstyle Y}  \quad (\alpha^*_{\scriptscriptstyle X} \beta^*_{\scriptscriptstyle Y}>0 \ \text{a.e.}) \iff \piXY,\piXYstar \in [\!\![\dep ]\!\!], \label{eqn:equiv2}\end{equation}
where we define $[\!\![\dep ]\!\!]$ as the equivalence class containing all the distributions sharing the dependence $\dep$. Accordingly, the dependence can be thought of as a maximal invariant under the group operation (\ref{eqn:equiv}). This allows the identification of a universal representation of $\dep$.

\subsection{A universal representation of $\dep$} \label{sec:universal}

Let $\Gamma_{\Xss \times \Yss}$ be the space of all finite signed measures $\gamma$ on $(\Xss \times \Yss,\Bs_\Xss \otimes \Bs_\Yss)$, and let 
\begin{equation} \Gamma_{\Ss_{XY}}^\circ = \{\gamma \in \Gamma_{\Xss \times \Yss}: |\gamma| \ll \piXY, \gamma(A,\Ss_Y) =  \gamma(\Ss_X,B) = 0\ \  \forall A \in \Bs_\Xss, \forall B \in \Bs_\Yss\}, \label{eqn:Gamma} \end{equation}
the closed linear subspace of $\Gamma_{\Xss \times \Yss}$ containing the measures absolutely continuous with respect to $\piXY$ and with marginals both identically null -- the notation $\Gamma_{\Ss_{XY}}^\circ$ makes it clear that this set is exclusively determined by the support $\Ss_{XY}$ of $\piXY$. Equipped with the total variation norm $\|\gamma\| = |\gamma|(\Xss \times \Yss)$, $\Gamma_{\Xss \times \Yss}$ and $\Gamma_{\Ss_{XY}}^\circ$ are Banach spaces \citep[Section 5.1]{Folland13}. Assume that we may find a basis for $\Gamma_{\Ss_{XY}}^\circ$, that is, a countable collection $\{\gamma_1,\gamma_2,\ldots \} \in \Gamma_{\Ss_{XY}}^\circ$, such that any $\gamma \in \Gamma_{\Ss_{XY}}^\circ$ admits the representation
\[\gamma = \sum_{k} a_k \gamma_k \]
for a (unique) real vector $(a_1,a_2,\ldots)$ \citep[Chapter 1]{Albiac16}. If so, call the number of elements in $\{\gamma_1,\gamma_2,\ldots \}$ the dimension of $\Gamma_{\Ss_{XY}}^\circ$, denoted $\text{dim}(\Gamma_{\Ss_{XY}}^\circ)$ (which may be infinite). Identifying such a basis is relatively straightforward in cases of interest -- see Section \ref{subsec:partcases}.  Then we have:

\begin{theorem} \label{thm:maxinv} Let $\piXY \in \overline{\powerset}_{\Xss \times \Yss}$ be a distribution with density $\phiXY = \diff \piXY /\diff (\mu_\Xss \times \mu_\Yss)$ on its support $\Ss_{XY}$. For any basis $\{\gamma_1,\gamma_2,\ldots \}$ of the set $\Gamma_{\Ss_{XY}}^\circ$ (\ref{eqn:Gamma}), the collection 
	\begin{equation} \DDelta_{XY} \doteq \left(\iint \log \phiXY \diff \gamma_1, \iint \log \phiXY \diff \gamma_2, \ldots \right) \in \R^{\text{dim}(\Gamma_{\Ss_{XY}}^\circ)} \label{eqn:coll} \end{equation}
	is a maximal invariant under the group operation (\ref{eqn:equiv}). 
\end{theorem}
%\begin{proof} See Appendix. \end{proof}

\noindent It directly follows:

\begin{corollary} \label{cor:dep} Under the conditions of Theorem \ref{thm:maxinv}, the pair $(\Ss_{XY},\DDelta_{XY})$ provides a complete representation of the dependence $\dep$ of the vector $(X,Y)$. 
\end{corollary}

\ppn The representation $(\Ss_{XY},\DDelta_{XY})$ is unique up to one-to-one correspondences. In particular, different bases $\{\gamma_1,\gamma_2,\ldots \}$ lead to different $\DDelta_{XY}$'s, but all are equivalent representations of $\dep$. As $\{\gamma_k\} \in \Gamma_{\Ss_{XY}}^\circ$, we have by definition $\iint \log \phiX \diff \gamma_k=\iint \log \phiY \diff \gamma_k =0 \ \forall k$ for any $\piXY \in \overline{\powerset}_{\Xss \times \Yss}$, so that
\begin{align*} \DDelta_{XY} = \left(\iint \log \phiXY \diff \gamma_1, \iint \log \phiXY \diff \gamma_2, \ldots \right) & = \left(\iint \log \phiYgX \diff \gamma_1, \iint \log \phiYgX \diff \gamma_2, \ldots \right) \\ & = \left(\iint \log \phiXgY \diff \gamma_1, \iint \log \phiXgY \diff \gamma_2, \ldots \right).\end{align*}
Also, $\Ss_{XY} = \{(x,y) \in \Ss_X  \times \Yss:   \phiYgX(y|x) > 0\} = \{(x,y) \in \Xss  \times \Ss_Y:   \phiXgY(x|y) > 0\}$ (up to sets of measure 0), so the support $\Ss_{XY}$ can be written in terms of $\{\phiYgX\}$ or $\{\phiXgY\}$ alone. Hence we can identify the pair $(\Ss_{XY},\DDelta_{XY})$ to the mapping $\Delta$ introduced in Section \ref{sec:compatibility}:
\[(\Ss_{XY},\DDelta_{XY}) = \Delta(\phiYgX) = \Delta(\phiXgY), \]
as prescribed by (\ref{eqn:deltas}). In other words, examining only one of the sets $\{\phiYgX\}$ or $\{\phiXgY\}$ of a distribution $\piXY$ is sufficient for identifying the dependence $\dep = (\Ss_{XY},\DDelta_{XY})$ of $\piXY$, as per Fact \ref{fact:conddens}.

\ppn It is also clear that, if $\phiXY$ satisfies (\ref{eqn:indep2}), then all elements of $\DDelta_{XY}$ are null. Hence, if $X \indep Y$ then $\DDelta_{XY} \equiv \zero$ \uline{and} $\Ss_{XY} = \Ss_X \times \Ss_Y$ (rectangular support, Section \ref{subsec:regio}). Conversely, if $\DDelta_{XY} = \zero$ {and} $\Ss_{XY} = \Ss_X \times \Ss_Y$, then $X \indep Y$. However, we may have $\Ss_{XY} \neq \Ss_X \times \Ss_Y$ and ${\bm \Delta}_{XY} = \zero$: then $X$ and $Y$ are {not} independent.  This is called `{quasi-independence}' \citep{Goodman68}: 

\begin{definition} \label{def:quasi} Let $(X,Y) \sim \piXY$ with density $\phiXY$ on the support $\Ss_{XY}$. Then, $X$ and $Y$ are quasi-independent if there exist (a.e.-)positive functions $q_{\scriptscriptstyle X} \in L_1(\Xss,\mu_\Xss)$ and $q_{\scriptscriptstyle Y} \in L_1(\Yss,\mu_\Yss)$ such that $\phiXY(x,y) = q_{\scriptscriptstyle X}(x) q_{\scriptscriptstyle Y}(y) \indic{(x,y) \in \Ss_{XY}}$.
\end{definition}

\noindent Independence implies quasi-independence, while quasi-independence implies independence only on a rectangular support. If the support is non-rectangular, the distribution $\piXY$ shows quasi-independence when the source of non-independence is purely `regional', i.e.\ exclusively due to the shape of the support (Section \ref{subsec:regio}). Extending the results of \cite{Caussinus65} in discrete settings (in particular, his `{\it Th\'eor\`eme IX}' and the following `{\it Remarque}'), we conclude:
\begin{corollary} \label{cor:indep} Let $(X,Y) \sim \piXY$. Under the conditions of Theorem \ref{thm:maxinv}:
	\begin{enumerate}
		\item $X$ and $Y$ are quasi-independent $\iff \DDelta_{XY} = \zero$;
		\item $X$ and $Y$ are independent $\iff \DDelta_{XY} = \zero$ and $\Ss_{XY} = \Ss_X \times \Ss_Y$.
	\end{enumerate}
\end{corollary}

\noindent Finally, it seems important to emphasise the following obvious consequence of Theorem \ref{thm:equivclass}:

\begin{corollary} \label{cor:excl} Any valid representation or description of the dependence of a bivariate distribution must be maintained over any equivalence class of dependence $[\!\![\dep]\!\!]$ defined by (\ref{eqn:equiv2}). In other words, any valid representation of $\dep$ must be invariant under the group operation (\ref{eqn:equiv}) or, yet equivalently, it must be recoverable from $(\Ss_{XY},\DDelta_{XY})$ only.
\end{corollary}

\subsection{Particular cases} \label{subsec:partcases}

The representation of dependence provided by Corollary \ref{cor:dep}  aligns with known results in familiar situations. For less well-studied cases, it allows unequivocal identification of how to capture and/or describe the dependence. Below some of those particular cases are explored.

\subsubsection{Two discrete variables} \label{subsec:discr}

Let $\Xss = \Yss = \N_0 = \{0,1,2,\ldots\}$ and $\mu_\Xss = \mu_\Yss$ be the counting measure on $(\N_0,2^{\N_0})$, meaning that $\phiXY(x,y) = \P(X=x,Y=y) \doteq p_{xy}$, for $(x,y) \in \N_0 \times \N_0$. In this case, the set $\Gamma_{\Xss \times \Yss}$ introduced in Section \ref{sec:universal} is isometrically isomorphic to $\ell_1(\N_0 \times \N_0)$, which makes it easy to identify bases for $\Gamma_{\Ss_{XY}}^\circ$. 

\ppn \uline{$(a)$ $X$ and $Y$ have finite supports $\Ss_X$ and $\Ss_Y$, $\Ss_{XY} = \Ss_X \times \Ss_Y$.} Assume $\Ss_X = \{0,1,\ldots,R-1 \}$ and $\Ss_Y = \{0,1,\ldots,S-1 \}$, for $2 \leq R, S < \infty$, and $p_{xy} >0$ $\forall (x,y) \in  \{0,1,\ldots,R-1 \} \times \{0,1,\ldots,S-1 \} = \Ss_{XY}$. Then $\Gamma_{\Ss_{XY}}^\circ$ can be identified to the space $\Ms_{R \times S}^\circ$ of $(R \times S)$-matrices with all rows and columns summing to 0. It is algebraic routine to show that $\text{dim}(\Gamma_{\Ss_{XY}}^\circ)= \text{dim}(\Ms_{R \times S}^\circ) = (R-1)(S-1)$, so the dependence in $(X,Y)$ can be entirely described by this number of parameters -- the $\chi^2$-test of independence in such a table articulates around $(R-1)(S-1)$ degrees of freedom, indeed. With $e_{xy}$ the $(R \times S)$-matrix whose all entries are zero except the $(x,y)$th entry which is 1, one can show that $\{E^{(00)}_{xy} \doteq (e_{xy} - e_{x0} - e_{0y} + e_{00}),x=1,\ldots,R-1,y=1,\ldots,S-1\}$ is a basis of $\Ms_{R \times S}^\circ$. Thus, the $(R-1)(S-1)$ components of $\DDelta_{\scriptscriptstyle XY}$ in (\ref{eqn:coll}) can be 
\begin{equation} \log \frac{p_{00}p_{xy}}{p_{x0}p_{0y}} , \quad x=1,\ldots,R-1, y=1,\ldots,S-1, \label{eqn:OR0} \end{equation}
in which we recognise the collection of log-odds ratios with respect to the pivot point $(0,0)$ -- any other pivot point may be chosen \citep[equation (2.11)]{Agresti13}. Another basis of $\Ms_{R \times S}^\circ$ is $\{E^{(\text{loc})}_{xy} \doteq (e_{x-1,y-1} - e_{x,y-1} - e_{x-1,y} + e_{xy}),x=1,\ldots,R-1,y=1,\ldots,S-1\}$, which would describe the dependence in terms of the `local odds ratios' \citep[equation (2.10)]{Agresti13}:
\begin{equation} \log \omega_{xy} \doteq \log \frac{p_{x-1\,y-1}p_{xy}}{p_{x-1\,y}p_{x\,y-1}}, \quad x=1,\ldots,R-1, y=1,\ldots,S-1. \label{eqn:OR1}  \end{equation}
It is known that the two sets of odds-ratios are equivalent (one may be recovered from the other) -- and so would be any other set built from yet a different basis of $\Ms_{R \times S}^\circ$. In agreement with Corollary \ref{cor:indep}, any collection of such log-odds-ratios is identically null when $X \indep Y$. In the case $R=S=2$, the dependence is fully captured by any one-to-one function of 
\[ \log \omega = \log \omega_{11} = \log \frac{p_{00}p_{11}}{p_{10}p_{01}},\]
as announced in Example \ref{ex:bivBern}.

\ppn \uline{$(b)$ $X$ and $Y$ have finite supports $\Ss_X$ and $\Ss_Y$, $\Ss_{XY} \neq \Ss_X \times \Ss_Y$.} A non-rectangular support for a bivariate discrete distribution (`incomplete table') is known to cause specific issues \citep[Chapter 5]{Bishop75}. Here, the layout of the `forbidden' entries is simply incorporated in $\Gamma_{\Ss_{XY}}^\circ$. As an example, consider a $(3 \times 3)$-distribution $\piXY$ supported on 
\begin{equation} \Ss_{XY} \sim \begin{pmatrix} \times & 0 & 0 \\ 0 & \times & \times \\ 0 & \times & \times \end{pmatrix} \label{eqn:suppex} \end{equation}
($\times$'s stand for non-zero entries). It can be seen that the space $\Gamma_{\Ss_{XY}}^\circ$ is one-dimensional, with basis vector $E^{(\text{loc})}_{22}=e_{11}-e_{21}-e_{12}+e_{22}$. Hence, given the support (\ref{eqn:suppex}), the dependence in such a distribution is fully characterised by (any one-to-one transformation of) 
\begin{equation} \log \frac{p_{11}p_{22}}{p_{21}p_{12}}. \label{eqn:OR2}  \end{equation}
The value of $p_{00}$ does not enter the description of dependence; indeed, on support (\ref{eqn:suppex}) it is entirely fixed by the marginal distributions: $p_{00} = \P(X=0) = \P(Y=0)$ -- the regional dependence structure reflected by (\ref{eqn:suppex}) imposes these two probabilities to be equal.

\ppn As another example, consider $\piXY$ supported on 
\begin{equation} \Ss_{XY} \sim \begin{pmatrix} 0 & \times & \times \\ \times & 0 & \times \\ \times & \times & 0 \end{pmatrix}. \label{eqn:suppdiag} \end{equation}
As $\Ss_{XY}$ does not comprise any positive rectangular sub-table,  a representation of dependence in terms of (log-)odds-ratios such as (\ref{eqn:OR0}), (\ref{eqn:OR1}) or (\ref{eqn:OR2}) is not possible here. Though, Corollary \ref{cor:dep} applies: the space $\Gamma_{\Ss_{XY}}^\circ$ is again unidimensional, with basis vector  $e_{01}-e_{02}-e_{10}+e_{12}+e_{20}-e_{21}$. Thus, on support (\ref{eqn:suppdiag}), the dependence is entirely captured by the single parameter
\begin{equation}  \log \frac{p_{01}p_{12}p_{20}}{p_{02}p_{10}p_{21}}. \label{eqn:OR3} \end{equation}
In particular,  $X$ and $Y$ are quasi-independent if and only if (\ref{eqn:OR3}) is 0 (Corollary \ref{cor:indep}), that is, $p_{01}p_{12}p_{20} = p_{02}p_{10}p_{21}$, confirming \citet[equation (2.23)]{Goodman68}.

\ppn When the structural zeros are too prominent, the space $\Gamma_{\Ss_{XY}}^\circ$ may contain only the null element. Then, the dependence is entirely determined by the specific support. An obvious example is a diagonal support such as
\begin{equation} \Ss_{XY} \sim \begin{pmatrix} \times  & 0 & 0 \\ 0 & \times & 0 \\ 0 & 0 & \times \end{pmatrix}. \label{eqn:suppdiag2} \end{equation}
Clearly $\Gamma_{\Ss_{XY}}^\circ = \{\zero\}$, the $(3 \times 3)$-null matrix. In other words, the fact that $\piXY$ is supported on the main diagonal is what characterises the dependence entirely -- concretely, it imposes $X=Y$ -- without any need for refinement: $p_{00},p_{11}$ and $p_{22}$ are fixed by the margins, which need to be identical in this case.

\ppn \uline{$(c)$ $X$ and $Y$ have infinite (countable) supports $\Ss_X$ and $\Ss_Y$.} Assume $\Ss_X = \Ss_Y = \N_0$. If $\Ss_{XY} = \Ss_X \times \Ss_Y$ (rectangular support), then the dependence in $(X,Y)$ is entirely described by an (infinite countable) collection of log-odds ratios, such as (\ref{eqn:OR0}) or (\ref{eqn:OR1}). If $\Ss_{XY} \neq \N_0 \times \N_0$, then the layout of structural zeros is part of the dependence. In case of scattered zeros not forming any particular pattern, the dependence may well remain captured by an (infinite countable) collection of $(2 \times 2)$-odds ratios, which essentially `avoid' the structural zeros. More specific patterns may give rise to other peculiar descriptions of the dependence structure, following the basis of $\Gamma_{\Ss_{XY}}^\circ$ which we can form. 

\subsubsection{Two (absolutely) continuous variables} \label{subsec:cont}

Let $\Xss = \Yss = [0,1]$ and $\mu_\Xss = \mu_\Yss$ be the Lebesgue measure. Denote $f_{\scriptscriptstyle XY}$ the joint probability density of $(X,Y)$, which is essentially bounded for $\piXY \in \overline{\powerset}_{\Xss \times \Yss}$ (see Appendix \ref{app:A}). The set $\Gamma_{\Ss_{XY}}^\circ$ can now be identified with the subset of functions $\psi$ of $L_1([0,1]^2)$ such that $\psi(x,y) = 0 $ for $(x,y) \notin \Ss_{XY}$, $\int_{\Ss_X} \psi(x,y)\, \diff x = 0$ $\forall y \in \Ss_Y$ and $ \int_{\Ss_Y} \psi(x,y) \, \diff y = 0$ $\forall x \in \Ss_X$. The tensorised Haar system is a (Schauder) basis for $L_1([0,1]^2)$ \citep[Chapter 6]{Albiac16}. Appropriate amendments of it provide bases for $\Gamma_{\Ss_{XY}}^\circ$.

\ppn \uline{$(a)$ $\Ss_{XY} =[0,1]^2 $.} Call $\psi_{0}(x,y) = (\indic{x \in [0,1/2)}- \indic{x \in [1/2,1)})(\indic{y \in [0,1/2)}- \indic{y \in [1/2,1)})$, and define $\psi_{k ,{\bm \ell}}(x,y) = 2^k\psi_0(2^k (x,y) - {\bm \ell})$ for $k \in \N_0$ and ${\bm \ell} \in \{0,\ldots,2^k-1\}^2$. Then $\left\{\psi_{k ,{\bm \ell}} \right\}_{k,{\bm \ell}}$ is a basis of $\Gamma_{\Ss_{XY}}^\circ$ -- see that it is the usual bivariate tensorised Haar basis, amputated from the scaling function (the `father wavelet') in both $x$- and $y$-directions so as to comply with $\int_0^1 \psi(x,y)\, \diff x = 0$ and $ \int_{0}^1 \psi(x,y) \, \diff y = 0$.

\ppn Then the dependence in $(X,Y)$ is  characterised by the collection $\DDelta_{\scriptscriptstyle XY} = \left\{\iint \log f_{\scriptscriptstyle XY} \ \psi_{k, {\bm \ell}} \,\diff x \diff y\right\}_{k,{\bm \ell}}$ (Corollary \ref{cor:dep}). The support of $\psi_{k ,{\bm \ell}}$ is the square $2^{-k}({\bm \ell} + [0,1]^2)$, which it partitions into 4 sub-squares, taking the value $2^k$ on the top-right and bottom-left sub-squares and the value $-2^k$ on the top-left and bottom-right sub-squares. Thus: 
\begin{align} \iint \log f_{\scriptscriptstyle XY} \ & \psi_{k, {\bm \ell}} \,\diff x \diff y =  2^k\iint_{2^{-k}({\bm \ell} + [0,1/2)^2)} \log f_{\scriptscriptstyle XY}\,\diff x \diff y + 2^k \iint_{2^{-k}({\bm \ell} + [1/2,1)^2)} \log f_{\scriptscriptstyle XY}\,\diff x \diff y \notag \\ & - 2^k\iint_{2^{-k}({\bm \ell} + [0,1/2)\times [1/2,1))} \log f_{\scriptscriptstyle XY}\,\diff x \diff y - 2^k\iint_{2^{-k}({\bm \ell} + [1/2,1) \times [0,1/2))} \log f_{\scriptscriptstyle XY} \,\diff x \diff y , \label{eqn:HaarlogOR}\end{align}
for $k \in \N_0$ and ${\bm \ell} \in \{0,\ldots,2^k-1\}^2$. We may recognise in (\ref{eqn:HaarlogOR}) some sort of log-odds ratio reminiscent of (\ref{eqn:OR1}): here the 4 relevant `probabilities' are akin to {geometric means} of $f_{\scriptscriptstyle XY}$ (i.e., average on the log-scale) over the corresponding sub-squares. As these get from coarser to finer, the complete collection $\left\{\iint \log f_{\scriptscriptstyle XY} \ \psi_{k, {\bm \ell}} \,\diff x \diff y\right\}_{k,{\bm \ell}}$ provides an appealing `multi-scale' representation of dependence.

\begin{remark} Such `multi-scale' approach, through artificially created $(2 \times 2)$-tables partitioning $\Ss_{XY}$, resonates with recent ideas of \cite{Ma19,Zhang19} and \cite{Gorsky22} for independence testing. Yet, those papers consider test statistics involving regular probabilities of $(X,Y)$ belonging to any of the above sub-squares; e.g., in our notation, odds ratios like
	\[\frac{\iint_{2^{-k}({\bm \ell} + [0,1/2)^2)} f_{\scriptscriptstyle XY}\,\diff x \diff y \times \iint_{2^{-k}({\bm \ell} + [1/2,1)^2)}  f_{\scriptscriptstyle XY}\,\diff x \diff y}{\iint_{2^{-k}({\bm \ell} + [0,1/2)\times [1/2,1))}  f_{\scriptscriptstyle XY}\,\diff x \diff y \times \iint_{2^{-k}({\bm \ell} + [1/2,1) \times [0,1/2))}  f_{\scriptscriptstyle XY} \,\diff x \diff y}; \]
	e.g., see \citet[Definition 3.6]{Zhang19}. As opposed to (\ref{eqn:HaarlogOR}), though, these quantities are not invariant under the group transformation (\ref{eqn:equiv}), hence do not describe $\dep$ (Corollary \ref{cor:excl}). \qed
\end{remark}

\noindent It so appears that $\DDelta_{\scriptscriptstyle XY} = \left\{\iint \log f_{\scriptscriptstyle XY} \ \psi_{k, {\bm \ell}} \,\diff x \diff y\right\}_{k,{\bm \ell}}$ consists of the coefficients representing 
\begin{equation} \overline{\lambda}_{\scriptscriptstyle XY} \doteq \log f_{\scriptscriptstyle XY} - \int \log f_{\scriptscriptstyle XY} \, \diff x - \int \log f_{\scriptscriptstyle XY} \, \diff y + \iint \log f_{\scriptscriptstyle XY} \, \diff x \diff y \label{eqn:lambdaXY} \end{equation}
in the regular bivariate tensorised Haar basis. The function $\overline{\lambda}_{\scriptscriptstyle XY}$, called `{centred log-odds ratio function}' in \cite{Osius04,Osius09}, is thus in one-to-one correspondence with $\DDelta_{XY}$, and is an equally valid representation of the dependence in $(X,Y)$. Clearly $\overline{\lambda}_{\scriptscriptstyle XY}$ is akin to the association parameters in log-linear parametrisations of contingency tables \citep[Chapter 7]{Agresti13}. For any arbitrary choice of reference point $(x_0,y_0) \in [0,1]^2$, we can also define the `{odds ratio function}' \citep{Chen07,Chen21,Chen15}
\begin{equation} \Omega_{\scriptscriptstyle XY}(x,y;x_0,y_0) \doteq  \frac{f_{\scriptscriptstyle XY}(x,y)f_{\scriptscriptstyle XY}(x_0,y_0)}{f_{\scriptscriptstyle XY}(x,y_0)f_{\scriptscriptstyle XY}(x_0,y)}, \label{eqn:ORfunct} \end{equation}
which is in one-to-one correspondence with $\overline{\lambda}_{\scriptscriptstyle XY}$:
\begin{multline*}  \overline{\lambda}_{\scriptscriptstyle XY}(x,y)  = \log \Omega_{\scriptscriptstyle XY}(x,y;x_0,y_0) -\int \log \Omega_{\scriptscriptstyle XY}(x,y;x_0,y_0) \diff x \\ - \int \log \Omega_{\scriptscriptstyle XY}(x,y;x_0,y_0) \, \diff y + \iint \log \Omega_{\scriptscriptstyle XY}(x,y;x_0,y_0) \, \diff x \diff y,\end{multline*}
\[ \Omega_{\scriptscriptstyle XY}(x,y;x_0,y_0) = \exp\left(\overline{\lambda}_{\scriptscriptstyle XY}(x,y) - \overline{\lambda}_{\scriptscriptstyle XY}(x,y_0) - \overline{\lambda}_{\scriptscriptstyle XY}(x_0,y) + \overline{\lambda}_{\scriptscriptstyle XY}(x_0,y_0)\right).\]
Hence for any $(x_0,y_0) \in [0,1]^2$, $\Omega_{\scriptscriptstyle XY}(\cdot,\cdot;x_0,y_0)$ is also a valid representation of dependence. 

\ppn Furthermore, if we let $k \to \infty$ in (\ref{eqn:HaarlogOR}) -- i.e., if we increase the resolution indefinitely -- the relevant squares become infinitesimally small and the limiting log-odds ratios turn into
\[\left(\log f_{\scriptscriptstyle XY}(x,y) + \log f_{\scriptscriptstyle XY}(x+\diff x,y + \diff y) - \log f_{\scriptscriptstyle XY}(x,y+\diff y) - \log f_{\scriptscriptstyle XY}(x+\diff x,y) \right)/(\diff x \diff y).\]
If $f_{\scriptscriptstyle XY}$ is mixed-differentiable, then this leads to the `{local dependence function}' \citep{Holland87}, viz.
\begin{equation} \gamma_{\scriptscriptstyle XY}(x,y) \doteq \frac{\partial^2 \log f_{\scriptscriptstyle XY}(x,y)}{\partial x \partial y}, \qquad (x,y) \in [0,1]^2. \label{eqn:ldf} \end{equation}
It can be verified that $$\gamma_{\scriptscriptstyle XY}(x,y)  = \frac{\partial^2 \log \Omega_{\scriptscriptstyle XY}(x,y;x_0,y_0)}{\partial x \partial y} \ =  \frac{\partial^2 \overline{\lambda}_{\scriptscriptstyle XY}(x,y)}{\partial x \partial y},$$ so $\gamma_{\scriptscriptstyle XY}$ is in one-to-one correspondence with $\Omega_{\scriptscriptstyle XY}(\cdot,\cdot;x_0,y_0)$ and $\overline{\lambda}_{\scriptscriptstyle XY}$ and is, therefore, a valid representation of dependence as well. Indeed, for rectangular supports, \citet[Lemma 3.2]{Wang93} proved that $\gamma_{\scriptscriptstyle XY}$ is maximal invariant under `marginal replacements', which we showed to be equivalent to (\ref{eqn:equiv}). This local dependence function was investigated further in \cite{Jones96,Jones98,Jones03,Molenberghs97}, and it was shown that $\gamma_{\scriptscriptstyle XY} \equiv 0$ if and only if $X\indep Y$ \citep[Lemma 4.2]{Holland87}. Clearly, $\gamma_{\scriptscriptstyle XY}$ is the continuous analogue of the family of local (log-)odds ratios (\ref{eqn:OR1}) while $\Omega_{\scriptscriptstyle XY}$ corresponds to the odds ratios (\ref{eqn:OR0}) with respect to a pivot point. 

\ppn Consider now Sklar's representation (\ref{eqn:Sklar}) of the distribution $\piXY$. As the copula $C$ is itself the distribution of an absolutely continuous vector $(U,V)$, its dependence may entirely be described by its odds ratio function (\ref{eqn:ORfunct}) as well, viz.
\[ \Omega_{\scriptscriptstyle UV}(u,v;u_\circ,v_\circ) = \frac{c(u,v) c(u_\circ,v_\circ)}{c(u_\circ,v) c(u,v_\circ)},  \] 
where $(u_\circ,v_\circ) \in [0,1]^2$ is any arbitrary pivot point, and $c(u,v) \doteq \partial^2 C(u,v) / \partial u \partial v$ is the copula density. Differentiating (\ref{eqn:Sklar}) yields $\fXY(x,y) = c(\FX(x),\FY(y)) \fX(x) \fY(y)$, which, plugged into (\ref{eqn:ORfunct}), immediately gives
\begin{equation} \Omega_{\scriptscriptstyle XY}(x,y;x_\circ,y_\circ) = \Omega_{\scriptscriptstyle UV}(\FX(x),\FY(y);\FX(x_\circ),\FY(y_\circ)). \label{eqn:margOmega} \end{equation}
Thus $\Omega_{\scriptscriptstyle XY}$ and $\Omega_{\scriptscriptstyle UV}$ are identical up to the usual one-to-one marginal re-labelling $u = \FX(x)$ and $v=\FY(y)$, meaning that $\FXY$ and $C$ share the same $\dep$ in agreement with Section \ref{subsec:equiv}: the change in marginal scales does not affect the inner structure. As the marginals of a copula are fixed by definition ($C \in \Fs(\Us_{[0,1]},\Us_{[0,1]})$), it follows from Definition \ref{def:dep} that there is one and only one copula with a given dependence structure/odds ratio function $\Omega_{\scriptscriptstyle UV}$. In view of this one-to-one correspondence,  positing the copula $C$ in (\ref{eqn:Sklar}) is perfectly equivalent to positing the odds ratio function $\Omega_{\scriptscriptstyle XY}$ of $\FXY$ as this one is $\Omega_{\scriptscriptstyle UV}$ after proper marginal re-labelling (\ref{eqn:margOmega}). Thus, the dependence in $(X,Y)$ is unequivocally described by $C$ as anticipated in Example \ref{ex:bivCop}. 

\ppn More specifically, the dependence structure of a bivariate Gaussian vector $(X,Y)$ (Example \ref{ex:Gauss}) is thus driven by that of $(U= \Upphi(\frac{X-\mu_{\scriptscriptstyle X}}{\sigma_{\scriptscriptstyle X}}),V= \Upphi(\frac{Y-\mu_{\scriptscriptstyle Y}}{\sigma_{\scriptscriptstyle Y}}))$, where $\Upphi$ denotes the standard normal cumulative distribution function (Probability Integral Transforms, PIT). The distribution of $(U,V)$ is then the Gaussian copula with parameter $\rho$, whose odds ratio function (\ref{eqn:ORfunct}) is (with $u_0 = 1/2$ and $v_0 = 1/2$): 
\begin{equation*} \Omega_{\scriptscriptstyle UV}(u,v;\frac{1}{2},\frac{1}{2}) = \frac{c_\rho(u,v) c_\rho(\frac{1}{2},\frac{1}{2})}{c_\rho(u,\frac{1}{2}) c_\rho(\frac{1}{2},v)} = \exp\left(\frac{\rho}{1-\rho^2} \Upphi^{-1}(u) \Upphi^{-1}(v)\right), \quad  (u,v) \in [0,1]^2. \label{eqn:ORUV}\end{equation*}
This parametric form of $\Omega_{\scriptscriptstyle UV}$ naturally follows from the Gaussian assumption, with its exact behaviour entirely dictated by $\rho$ (given that $\rho/(1-\rho^2)$ is a one-to-one function of $\rho$). Therefore, in a Gaussian vector, Pearson's correlation $\rho$ fully describes the dependence, as announced in Example \ref{ex:Gauss}. In particular, $X \indep Y \iff U \indep V \iff \log \Omega_{\scriptscriptstyle UV} \equiv 0 \iff \rho =0$.

\ppn We may want to write the odds ratio function on the original support of $(X,Y)$. Around $x_0 = \mu_{\scriptscriptstyle X}$ and $y_0 = \mu_{\scriptscriptstyle Y}$, we obtain by substituting the bivariate Gaussian density in (\ref{eqn:ORfunct}):
\[\Omega_{\scriptscriptstyle XY}(x,y;\mu_{\scriptscriptstyle X},\mu_{\scriptscriptstyle Y}) = \exp\left(\frac{\rho}{1-\rho^2} \left(\frac{x-\mu_{\scriptscriptstyle X}}{\sigma_{\scriptscriptstyle X}}\right)\left(\frac{y-\mu_{\scriptscriptstyle Y}}{\sigma_{\scriptscriptstyle Y}}\right)\right), \qquad (x,y) \in \R^2, \]
showing the same univocal effect of $\rho$. It may be surprising to see the marginal parameters $(\mu_{\scriptscriptstyle X},\sigma_{\scriptscriptstyle X})$ and $(\mu_{\scriptscriptstyle Y},\sigma_{\scriptscriptstyle Y})$ appear in this representation of dependence, though. Yet, this is just the principle of equivariance exposed in Section \ref{subsec:equiv} coming into play. In fact, for any two bimeasurable bijections $\Phi,\Psi$, we have (similarly to (\ref{eqn:margOmega})):
\begin{equation} \Omega_{\scriptscriptstyle \Phi(X)\Psi(Y)}(x',y';x'_\circ,y'_\circ) =  \Omega_{\scriptscriptstyle XY}(\Phi^{-1}(x'),\Psi^{-1}(y');\Phi^{-1}(x'_\circ),\Psi^{-1}(y'_\circ)). \label{eqn:omtrans} \end{equation}
Any arbitrary parameters may hide into $\Phi$ and $\Psi$, without affecting the dependence captured by $\Omega_{\scriptscriptstyle XY}$. With the above choice of PIT's, $\Phi$ and $\Psi$ involve the marginal parameters $(\mu_{\scriptscriptstyle X},\sigma_{\scriptscriptstyle X})$ and $(\mu_{\scriptscriptstyle Y},\sigma_{\scriptscriptstyle Y})$, but this should not be confused with some interference between these parameters and the dependence. 

\ppn Now, it is shown in \cite{Jones96} that $\gamma_{\scriptscriptstyle XY}(x,y)$ may be interpreted as an (appropriately scaled) localised version of Pearson's correlation around the point $(x,y)$ -- hence the name `{local dependence}'. In general, such local dependence is not uniform over the whole support $\Ss_{XY}$; e.g., some distributions show `tail dependence' \citep[Section 2.13]{Joe2015}, meaning that the variables tend to be more strongly tied in their areas of extreme values than around their median point. Yet, for a Gaussian vector, we see that   
\begin{equation*} \gamma_{\scriptscriptstyle XY}(x,y)  = \frac{\partial^2 \log \Omega_{\scriptscriptstyle XY}(x,y;\mu_{\scriptscriptstyle X},\mu_{\scriptscriptstyle Y})}{\partial x \partial y} = \frac{\rho}{1-\rho^2} \frac{1}{\sigma_{\scriptscriptstyle X} \sigma_{\scriptscriptstyle Y}} \qquad \forall (x,y) \in \R^2,   
	\label{eqn:ldfGauss} \end{equation*}
indicating that the local dependence {\it is} uniform over the whole of $\R^2$. In agreement with (\ref{eqn:decomp}), we can mount arbitrary $\R$-supported marginals on this specific dependence structure for creating a multitude of distributions having this trait in common with the bivariate Gaussian. 

\ppn \uline{$(b)$ $\Ss_{XY} \subsetneq  [0,1]^2$ (non-rectangular support).} The support $\Ss_{XY}$ of any $\piXY \in \overline{\powerset}_{\Xss \times \Yss}$ can be expressed as a countable union of rectangles (see Appendix \ref{app:A}). Thus, at resolution high enough, the elements of a Haar basis for $\Gamma_{\Ss_{XY}}^\circ$ may cover the whole of $\Ss_{XY}$ arbitrarily close by a rectangular pavement. In particular, if we push the resolution to the limit $k \to \infty$, then (assuming $f_{\scriptscriptstyle XY}$ is mixed-differentiable) we recover the local dependence function (\ref{eqn:ldf}) on $\Ss_{XY}$, irrespective of the shape of that set. In agreement with Corollary \ref{cor:dep}, the couple $(\Ss_{XY},\gamma_{\scriptscriptstyle XY})$ provides a complete characterisation of the dependence in $(X,Y)$. E.g., the vector $(X,Y)$ will show quasi-independence on its support $\Ss_{XY}$ if and only if $\gamma_{\scriptscriptstyle XY} \equiv 0$ on $\Ss_{XY}$ (Corollary \ref{cor:indep}).  

\subsubsection{One (absolutely) continuous and one discrete variable}

Let  $\Xss = \N_0 = \{0,1,\ldots\}$, $\mu_\Xss$ be the counting measure, $\Yss = [0,1]$ and $\mu_\Yss$ be the Lebesgue measure. Assume for simplicity that $\piXY$ is supported on $\Ss_{XY} = \N_0 \times [0,1]$ (rectangular support). In such a {mixed} vector, the density $\phiXY$ takes the form
\[\phiXY(x,y) = \P(X=x|Y=y) f_{\scriptscriptstyle Y}(y) = f_{\scriptscriptstyle Y|X}(y|x) \P(X=x), \qquad (x,y) \in \N_0 \times [0,1], \]
where the notations $\P(X=\cdot|Y=\cdot)$, $f_{\scriptscriptstyle Y|X}(\cdot|\cdot)$, $\P(X=\cdot)$ and $f_{\scriptscriptstyle Y}(\cdot)$ are self-evident. Ideas from Sections \ref{subsec:discr} and \ref{subsec:cont} can be combined to identify a basis for $\Gamma_{\Ss_{XY}}^\circ$, of the form $\{ (\delta_x - \delta_{x-1})\psi_{k \ell}: x=1,2,\ldots,k \in \N_0, \ell \in \{0,\ldots,2^k-1\}\}$, where $\delta_x$ is the unit mass at $x \in \N_0$, and for $k \in \N_0, \ell \in  \{0,\ldots,2^k-1\}$, $\psi_{k \ell}(y) =2^{k} \psi_0(2^{k} y - \ell)$ with $\psi_0(y) = \indic{y \in [0,1/2)}- \indic{y \in [1/2,1)}$. The dependence in $(X,Y)$ is then characterised (Corollary \ref{cor:dep}) by the collection 
\begin{multline} \left\{2^{k} \int_{{2^{-k}(\ell + [0,1/2))}} \log \frac{\P(X=x|Y=y)}{\P(X=x-1|Y=y)} \, \diff y  - 2^{k} \int_{{2^{-k}(\ell + [1/2,1))}} \log \frac{\P(X=x|Y=y)}{\P(X=x-1|Y=y)} \, \diff y\right\}, \label{eqn:discrcontcoll} \end{multline}
for $x=1,2,\ldots,k \in \N_0, \ell \in \{0,\ldots,2^k-1\}$. If we increase arbitrarily the resolution level of the univariate Haar system $\{\psi_{k \ell}\}$ ($k \to \infty$), the elements of (\ref{eqn:discrcontcoll}) become
\begin{equation} \gamma_{\scriptscriptstyle XY}(x,y) = \frac{\partial  }{\partial y} \log \frac{\P(X=x|Y=y)}{\P(X=x-1|Y=y)}, \qquad (x,y) \in \{1,2\ldots\} \times [0,1], 
	\label{eqn:ldfmixed} \end{equation}
or equivalently
\[\gamma_{\scriptscriptstyle XY}(x,y) = \frac{\partial  }{\partial y} \log \frac{f_{\scriptscriptstyle Y|X}(y|x)}{f_{\scriptscriptstyle Y|X}(y|x-1)}, \qquad (x,y) \in \{1,2\ldots\} \times [0,1], \]
assuming that $\P(X=x|Y=\cdot)$ and $f_{\scriptscriptstyle Y|X}(\cdot|x)$ are differentiable as functions of $y$ for all $x \in \N_0$. In such a case, either of these two functions represents unequivocally the dependence in $(X,Y)$.

\begin{example} \label{ex:logistic} Consider the simple logistic regression framework \citep[Chapter 4]{Agresti13}, with $X \in \{0,1\}$ and $\P(X=1|Y=y) = \exp(\alpha + \beta y) / (1+
	\exp(\alpha + \beta y))$, for $\alpha,\beta \in \R$. With $X \in \{0,1\}$, (\ref{eqn:ldfmixed}) reduces down to a univariate function $\gamma_{\scriptscriptstyle XY}(y)$, which is here constant:
	\[\gamma_{\scriptscriptstyle XY}(y) =  \frac{\partial  }{\partial y} \log \exp(\alpha + \beta y) = \beta. \]
	Therefore, under this model, the dependence between $X$ and $Y$ is captured by one single parameter, which should be any one-to-one transformation of the slope parameter $\beta$. The constancy of $\gamma_{\scriptscriptstyle XY}$ in this case complements the results of \cite{Jones98}, who showed that, in continuous situations, the local dependence function (\ref{eqn:ldf}) is constant in the case of a generalised linear model (GLM) with canonical link -- this includes the bivariate Gaussian case of Section \ref{subsec:cont}. We note that, effectively, (\ref{eqn:ldfmixed}) can be regarded as a `mixed' version of (\ref{eqn:ldf}). 
	
	\ppn By contrast, in the probit model $\P(X=1|Y=y) = \Upphi(\alpha +\beta y)$, $\gamma_{\scriptscriptstyle XY}(y)$ is not constant and involves both $\alpha$ and $\beta$, making dependence in `logit' and `probit' vectors different. This could be inferred directly from Definition \ref{def:dep}: given that logit and probit links may be used for defining different joint distributions for $(X,Y)$ with the same marginals $X \sim \text{Bern}(p)$ and $Y \sim \piY$, the difference between both models must be distinct dependence structures. More generally, it appears that the choice of a link function in a GLM is essentially one of dependence structure between response and predictor, with the canonical link accounting for a constant local dependence entirely describable by a single parameter. \qed
\end{example}

\subsubsection{Hybrid variables}

By `hybrid' we mean a distribution which is a mixture of a discrete component and an absolutely continuous component. For illustration, consider the case of a vector $(X,Y)$ where both $\piX$ and $\piY$ are made up of a probability atom at 0 and an absolutely continuous part over $\R^+$. These marginal distributions are absolutely continuous with respect to the reference measures $\mu_\Xss = \mu_\Yss = \delta_0 + \lambda$ (where $\delta_0$ is the Dirac measure at 0 and $\lambda$ is the Lebesgue measure on $\R^+$). The density $\phiXY$ of such a vector takes the form
\[\phiXY(x,y) = \left\{ \begin{array}{l l} p_{00} & \text{ if } x=0,y=0 \\ p_{10} f_{\scriptscriptstyle X0}(x) & \text{ if } x > 0, y = 0 \\ p_{01} f_{\scriptscriptstyle 0Y}(y) & \text{ if } x = 0, y > 0 \\ p_{11} f_{\scriptscriptstyle XY}(x,y) & \text{ if } x > 0, y > 0 \end{array}\right. ,\]
where $p_{00}, p_{10}, p_{01}, p_{11}$ are probabilities summing to 1, $f_{\scriptscriptstyle X0}$ is the conditional density of $X$ given $(X> 0, Y=0)$, $f_{\scriptscriptstyle 0Y}$ is the conditional density of $Y$ given $(X = 0, Y > 0)$ and $f_{\scriptscriptstyle XY}$ is the conditional density of $(X,Y)$ given that $(X> 0,Y> 0)$. We assume that $f_{\scriptscriptstyle X0}(x) >0, f_{\scriptscriptstyle 0Y}(y)>0$ and $f_{\scriptscriptstyle XY}(x,y) >0$ $\forall (x,y) \in \R^+ \times \R^+$ (rectangular support). In this case, the most natural representation of the dependence in $(X,Y)$ seems to be the odds-ratio function with respect to the pivot point $(0,0)$:
\[\Omega_{\scriptscriptstyle XY}(x,y;0,0) = \left\{ \begin{array}{l l} \frac{p_{00}p_{11}}{p_{10}p_{01}} \frac{f_{\scriptscriptstyle XY}(x,y)}{f_{\scriptscriptstyle X0}(x) f_{\scriptscriptstyle 0Y}(y)} & (x,y) \in \R^+ \times \R^+  \\ 1 & (x,y) \in \{x=0\} \cup \{y=0\} \end{array} \right. .\]
We recognise the usual odds ratio $(p_{00}p_{11})/(p_{10}p_{01})$ formed from the $(2 \times 2)$-table $(X=0,X > 0) \times (Y=0, Y > 0)$, as well as a continuous contribution for $(x,y) \in \R^+ \times \R^+$.

\section{Dependence versus concordance} \label{sec:concord}

Most of statistical theory developed from the initial works of Galton and Pearson on the Gaussian distribution \citep{Stigler02}. As a result, much of the statistical jargon in use today retains a strong Gaussian connotation. For example, the fact that the dependence reduces down to the parameter $\rho$ in a bivariate Gaussian distribution (Example \ref{ex:Gauss}) largely explains why the correlation coefficient has been referred to as a `dependence measure' universally, including outside the Gaussian model. In fact, a bivariate Gaussian vector $(X,Y)$ is so defined that any relation which may exist between $X$ and $Y$ must agree with $\rho$. This has caused, through a similar process of conflation, most of those relations to be indiscriminately labelled as `dependence'. Here we elaborate on the specific case of `concordance'.

\ppn Two variables $X$ and $Y$ show {concordance} when large (small) values of $X$ occur typically with large (small) values of $Y$; i.e., when $X$ and $Y$ vary in the same direction: increasing the value of one tends to increase the value of the other. {Discordance} is the reverse effect: $X$ and $Y$ vary in opposite directions (increasing one tends to decrease the other). If $X$ and $Y$ are concordant, then $X$ and $-Y$ are discordant; and discordance may be regarded as `negative concordance'. Note that, for introducing `small' and `large' values for $X$ and $Y$, we need to explicitly assume in this section that $\Xss$ and $\Yss$ are ordered sets; say, $\Xss$ and $\Yss$ are subsets of $\R$. 

\ppn Concordance is formally defined through a partial ordering of the distributions of a given Fr\'echet class. %Define the {joint cumulative distribution function} (cdf) of $\piXY$ as
%\[ F_{\scriptscriptstyle XY}(x,y) \doteq \pi_{\scriptscriptstyle XY}\big((-\infty,x] \times (-\infty,y] \big),\quad (x,y) \in \Xss \times \Yss \]
%and its marginal counterparts $F_{\scriptscriptstyle X}(x) = \piX\big((-\infty,x]\big)$ and $F_{\scriptscriptstyle Y}(y) = \piY\big((-\infty,y]\big)$. 
Specifically, for $\pi^{(1)}_{\scriptscriptstyle XY}$, $\pi^{(2)}_{\scriptscriptstyle XY} \in \Fs(\piX,\piY)$ with respective cdf's $F^{(1)}_{\scriptscriptstyle XY}$ and $F^{(2)}_{\scriptscriptstyle XY}$, we say that $\pi^{(1)}_{\scriptscriptstyle XY}$ is `{less concordant}' than  $\pi^{(2)}_{\scriptscriptstyle XY}$, denoted $\pi^{(1)}_{\scriptscriptstyle XY} \precsim_c \pi^{(2)}_{\scriptscriptstyle XY}$, when 
\begin{equation} F^{(1)}_{\scriptscriptstyle XY}(x,y) \leq F^{(2)}_{\scriptscriptstyle XY}(x,y), \quad \forall (x,y) \in \Xss \times \Yss \label{eqn:concord} \end{equation}
\citep{Scarsini84,Tchen80}. Any numerical description of bivariate distributions which happens to be monotonic with respect to the ordering $\precsim_c$ can, therefore, be thought of as a measure of concordance. This is the case of Spearman's $\rho_S$ and Kendall's $\tau$, explicitly defined as a difference between probabilities of `concordance' and `discordance' \citep{Hoeffding47,Kruskal58}, viz.
\begin{align*} \rho_S & =3\left(\P((X_1-X_2)(Y_1-Y_3) >0) - \P((X_1-X_2)(Y_1-Y_3) <0)\right),  \\ \tau & = \P((X_1-X_2)(Y_1-Y_2) >0) - \P((X_1-X_2)(Y_1-Y_2) <0),  \end{align*}
for independent copies $(X_1,Y_1),(X_2,Y_2),(X_3,Y_3) \sim \piXY$. Indeed $\rho_S$ and $\tau$  are such that: 
\begin{equation*} \pi^{(1)}_{\scriptscriptstyle XY} \precsim_c \pi^{(2)}_{\scriptscriptstyle XY} \quad \Rightarrow \quad \left\{\begin{array}{l} \rho_S^{(1)} \leq \rho_S^{(2)} \\  \tau^{(1)} \leq \tau^{(2)} \end{array} \right. \label{eqn:concordrhotau} \end{equation*} 
\citep[Corollary 3.2]{Tchen80}.  As in Section \ref{subsec:equiv}, introduce $(\Xss',\Bs_{\Xss'})$ and $(\Yss',\Bs_{\Yss'})$ two measurable spaces isomorphic to $(\Xss,\Bs_{\Xss})$ and $(\Yss,\Bs_{\Yss})$, and define two bimeasurable bijections $\Phi: \Xss \to \Xss'$ and $\Psi:\Yss \to \Yss'$. Then, for $k=1,2$, let $\pi^{(k)}_{\scriptscriptstyle \Phi(X)\Psi(Y)}(A' \times B') = \pi^{(k)}_{\scriptscriptstyle XY}(\Phi^{-1}(A') \times\Psi^{-1}(B'))$, $\forall (A',B') \in \Bs_{\Xss'} \otimes \Bs_{\Yss'}$. If $\Phi$ and $\Psi$ are order-preserving (that is, if they are increasing), then $\FXY(x,y) = F_{\scriptscriptstyle \Phi(X)\Psi(Y)}(\Phi(x),\Psi(y))$ $\forall (x,y) \in \Xss \times \Yss$, and it is clear from (\ref{eqn:concord}) that  
\begin{equation*} \pi^{(1)}_{\scriptscriptstyle \Phi(X)\Psi(Y)} \precsim_c \pi^{(2)}_{\scriptscriptstyle \Phi(X)\Psi(Y)} \quad \iff \quad \pi^{(1)}_{\scriptscriptstyle XY} \precsim_c \pi^{(2)}_{\scriptscriptstyle XY}. \label{eqn:concor} \end{equation*}
%Furthermore, $\rho_S$ and $\tau$ in (\ref{eqn:rhoS})-(\ref{eqn:tau}) may be equivalently defined as 
%\begin{align*} \rho_S & =3\left(\P((\Phi(X_1)-\Phi(X_2))(\Psi(Y_1)-\Psi(Y_3)) >0) - \P((\Phi(X_1)-\Phi(X_2))(\Psi(Y_1)-\Psi(Y_3)) <0)\right), \\ \tau & = \P((\Phi(X_1)-\Phi(X_2))(\Psi(Y_1)-\Psi(Y_2)) >0) - \P((\Phi(X_1)-\Phi(X_2))(\Psi(Y_1)-\Psi(Y_2)) <0). \end{align*}
In other words, concordance is a trait invariant to increasing transformations of the marginals. Evidently $\rho_S$ and $\tau$ are invariant under such transformations, too.

\ppn Now, the Gaussianity of $(X,Y)$ parameterised as in (\ref{eqn:paramGauss}) imposes the linear relationship 
\[Y = \mu_{\scriptscriptstyle Y} + \rho \frac{\sigma_{\scriptscriptstyle Y}}{\sigma_{\scriptscriptstyle X}} (X - \mu_{\scriptscriptstyle X}) + \epsilon, \qquad \text{ with } \epsilon \sim  \Ns(0,\sigma^2_{\scriptscriptstyle X}(1-\rho^2)).\]
Clearly, a positive (resp.\ negative) slope for this straight line implies concordance (resp.\ discordance) between $X$ and $Y$. Hence the value of $\rho$ (i.e., the dependence) directly dictates the concordance/discordance in $(X,Y)$ -- both through its sign and its absolute value, as a reduced variability for $\epsilon$ strengthens the sense of concordance/discordance between $X$ and $Y$. So, it is not surprising that $\rho_S$ and $\tau$ are one-to-one functions of $\rho$ in this model, and thus valid representations of $\dep$ (Example \ref{ex:Gauss}). %\footnote{On a side note, the only increasing marginal transformations $\Phi$ and $\Psi$ which are allowed {in order to stay in the Gaussian model} are linear transformations, to which $\rho$ is known to be invariant.} 
That the two concepts coincide in the Gaussian case may explain why concordance is commonly understood as `some form' of dependence. For example, the term `monotone dependence' (`positive'/`negative', for concordance/discordance) is often used, while Spearman's $\rho_S$ and Kendall's $\tau$ are routinely referred to as `{dependence measures}', in spite of specifically accounting for {concordance}. Yet, in general, dependence (understood through Definition \ref{def:dep}) and concordance may well be unrelated.

\ppn Concordance/discordance explicitly implies an `influence' of one variable on the other (`{if $X$ increases, $Y$ tends to increase/decrease}'), hence non-independence (Section \ref{subsec:conddensdep}). Accordingly, $X  \indep Y$ implies $X$ and $Y$ not in concordance/discordance. But $X$ and $Y$ not in concordance/discordance does not imply their independence \citep[Theorem 4]{Scarsini84}. Hence the same concordance status may exist with different dependence structures (e.g., `not concordant/discordant' may be observed in vectors showing independence as well as in vectors showing non-independence). Reversely, the same dependence can exist with different levels of concordance, as illustrated by the following example.

\begin{example} Let $X \sim  \Us_{[0,1]}$ and $Y = ZX$, where $Z$ is a Rademacher random variable\footnote{$\frac{Z+1}{2} \sim \text{Bern}(1/2)$.} independent of $X$. The distribution $\piXY$ is supported on the two diagonals of the unit square:\footnote{The argument is clearer when examining this singular distribution, although it may not perfectly fit in the general framework ($\piXY \not\in \powerset_{\Xss \times \Yss}$). An absolutely continuous version may easily be constructed, though.} call the 4 arms of this cross `top right', `bottom right', `bottom left' and `top left'. Due to the perfect symmetry, $(X,Y)$ cannot show any concordance or discordance, and indeed, $\tau = 0$ in this case. Now, define $\piXYstar$ through (\ref{eqn:equiv}) with $\alpha^*_{\scriptscriptstyle X}(x) = 2(1-a)\indic{x \in [0,1/2]}+2a\indic{x \in (1/2,1]}$ and $\beta^*_{\scriptscriptstyle Y}(y) = 2(1-b)\indic{y \in [0,1/2]}+2b\indic{y \in (1/2,1]}$, for $a,b \in (0,1)$. According to Theorem \ref{thm:equivclass}, this maintains the dependence ($\dep = \dep^*$) regardless of the exact values of $a$ and $b$ -- in particular, the support is the same. However, it can be verified that Kendall's tau is $\tau^* = (2a-1)(2b-1)$ for $\piXYstar$. If both $a$ and $b$ are close to 1, $\piXYstar$ is mostly concentrated on the top right arm of the cross, and thus shows strong concordance ($\tau^* \simeq 1$). By contrast, if $a$ is close to 1 but $b$ is close to 0, then $\piXYstar$ is mostly concentrated on the bottom right arm, and thus shows strong discordance ($\tau^* \simeq -1$). In fact, varying the values of $a,b \in (0,1)$, we can obtain any value of $\tau^* \in (-1,1)$ while keeping the initial dependence $\dep$. \qed
\end{example} 

\ppn Thus, beyond the triviality `independence' $\Rightarrow$ `no concordance/discordance', there is in general no link between concordance and dependence. This can be appreciated pragmatically: for defining $\dep$ we could remain oblivious of the ordering in $\Xss$ and $\Yss$ -- in fact, we allowed such an ordering to be entirely reshuffled in Section \ref{subsec:equiv} -- hence concordance, exclusively based on such ordering, must be of different nature. For example, the dependence in a discrete distribution is characterised by a collection of log-odds ratios such as (\ref{eqn:OR0}), which each captures how likely the specific values $X=x$ and $Y=y$ are to occur together. Clearly, these log-odd ratios can be listed in any order for covering all combinations of $X$- and $Y$-values in $\Ss_{XY}$. In fact, we could describe the dependence between two categorical {\it non-ordinal} variables by a similar collection of odds-ratios. The same observation may be made in the continuous case as well, from the odds ratio function $\Omega_{\scriptscriptstyle XY}$ (\ref{eqn:ORfunct}): the assumed order on $\Xss$ and $\Yss$ may be arbitrarily altered through one-to-one bimeasurable Lebesgue-measure preserving functions $\Phi: \Xss \to \Xss$ and $\Psi: \Yss \to \Yss$ while keeping the same information about the dependence {\it \`a la} (\ref{eqn:omtrans}). Notably, this agrees with \cite{Renyi59}'s views, which require valid dependence measures to be invariant under one-to-one (but not necessarily increasing) transformations of the margins (his axiom F).

\ppn Now, some specific dependence structures do dictate the concordance. For example, a vector $(X,Y) \sim \piXY$ is said to show `{Positive Likelihood Ratio Dependence}' (PLRD) \citep[Section 8]{Lehmann66} when, for all $x_1,x_2 \in \Xss$ and $y_1,y_2 \in \Yss$ such that $x_1<x_2$ and $y_1 < y_2$,
\begin{equation} \phiXY(x_1,y_1)\phiXY(x_2,y_2) \geq \phiXY(x_2,y_1)\phiXY(x_1,y_2)  \label{eqn:PLR} \end{equation}
(with strict inequality for at least one pair $(x_1,y_1),(x_2,y_2)$). Evidently (\ref{eqn:PLR}) remains true for the density $\phiXYstar$ of any distribution $\piXYstar$ satisfying (\ref{eqn:equiv}), hence PLRD is truly a dependence property: in a given equivalence class $[\!\![\dep ]\!\!]$, either all distributions show PLRD, or none does. It is clear that (\ref{eqn:PLR}) induces a sense of concordance, and indeed, PLRD implies what has been called `{Positive Quadrant Dependence}' (PQD) \citep[Theorem 2.3]{Joe97}, which itself implies $\rho_S>0$ and $\tau >0$ \citep[Corollary 1]{Lehmann66}. Note that PQD, defined as when $F_{\scriptscriptstyle XY}(x,y) \geq F_{\scriptscriptstyle X}(x) F_{\scriptscriptstyle Y}(y)$ $\forall (x,y) \in \Xss \times \Yss$, is {not} necessarily preserved by (\ref{eqn:equiv}), hence is {\it not} a dependence property: within the same equivalence class $[\!\![\dep ]\!\!]$, some distributions may show PQD while others don't, in violation of Corollary \ref{cor:excl}. The fact that all bivariate Gaussian distributions with $\rho \neq 0$ show either PLRD ($\rho >0$) or its negative counterpart ($\rho < 0$) \citep[Example 10]{Lehmann66} explains why dependence and concordance merge in that case. The same may be said about the bivariate Bernoulli model (Example \ref{ex:bivBern}), from \citet[Example 11]{Lehmann66}. 

%\ppn This suggests another definition of the copula $C$ of $\FXY$: this is the {\it unique} distribution with the same dependence structure and the same concordance as $\FXY$ but with uninformative uniform margins. The condition `same concordance' is what makes this distribution unique, as it is known that $\Phi = \FX$ is the only {\it increasing} function such that $\Phi(X) \sim \Us_{[0,1]}$; e.g., see \citet[Theorem 1.5.1]{Panaretos20}. Other one-to-one functions of $X$ and $Y$ may produce $\Us_{[0,1]}$-marginals (think of $1-\FX(X)$ and $1-\FY(Y)$, for example), without affecting the dependence structure of $(X,Y)$ \citep[Section 3.3]{Geenens22}. Yet, increasing transformations are the only ones which preserve concordance \citep[Section 5]{Geenens22}.

\section{Measuring the dependence} \label{sec:depmes}

\subsection{General comments}

It appears clearly from Corollary \ref{cor:dep} that dependence is a polymorphic concept that cannot be entirely characterised by one single parameter in general (parametrically restricted environments such as in Examples \ref{ex:Gauss}, \ref{ex:bivBern} and \ref{ex:logistic} are exceptions). This said, it often remains essential to express succinctly, by a suitable index, how strongly two variables `influence' each other. This is the main motivation behind the numerous {dependence measures} mentioned in Section \ref{sec:intro}. Yet, most of those measures do {not} actually deal with dependence as delineated by Definition \ref{def:dep}. For example, Section \ref{sec:concord} described how {concordance} measures such as Spearman's $\rho_S$ or Kendall's $\tau$ may hardly be thought of as measuring {dependence}, despite being non-null only in case of non-independence. In particular, $\rho_S = 0$ or $\tau = 0$ does not mean `no dependence'.

\ppn Failure to identify independence is not the only impediment, though. Consider the Mutual Information
\begin{equation} \textup{MI}(X,Y) = \Is(\piXY \| \piX \times \piY) = \iint \log \frac{\diff \piXY}{\diff (\piX \times \piY)} \diff \piXY,   \label{eqn:MI0}\end{equation}
a popular `dependence measure' owing to the fact that $\textup{MI}(X,Y) = 0  \iff X \indep Y$ \citep[p.\,28]{Cover06}. Yet, it appears directly from (\ref{eqn:MI0}) that $\textup{MI}(X,Y)$ is not invariant under the group operation (\ref{eqn:equiv}); meaning that it may vary from a distribution to another within the same equivalence class of dependence $[\!\![\dep]\!\!]$. Hence the Mutual Information violates Corollary \ref{cor:excl} and cannot be thought of as quantifying {dependence} as such. In fact, as its name suggests, $\textup{MI}(X,Y)$ measures the {\it information} shared between $X$ and $Y$, a specific notion defined by \citet{Shannon1948} and interpreted in terms of the reduction of the uncertainty surrounding $Y$ once we observe $X$ (or vice-versa). While such shared information may only be null when $X \indep Y$, this concept should not be amalgamated with dependence. For example, \citet[p.\ 643]{Geenens20} exhibit a vector $(X,Y)$ in which $X$ and $Y$ are `nearly-independent' while their Mutual Information is arbitrarily large, illustrating the discrepancy between the value of $\textup{MI}(X,Y)$ and the underlying dependence structure.

\ppn Similar observations may be made about many other alleged `dependence measures'. It can reasonably be admitted that, if $X$ and $Y$ are independent, then `nothing' may exist between them. Any index quantifying the magnitude of a trait which cannot be present between two independent variables is automatically null for $X \indep Y$; and consequently, a non-null value for that index reveals non-independence. Yet, this does not mean that the said index relates to Definition \ref{def:dep}. If it does not, calling it a `dependence measure' is misleading -- such an index should rather be called a `{\it non-independence indicator}'.\footnote{Note that it would be perfectly valid to use such indicator as the test statistic for independence testing.}

\ppn In fact, Corollary \ref{cor:excl} imposes that any valid description of the dependence in the vector $(X,Y)$ -- which includes any attempt to quantify its strength -- must be based on $(\Ss_{XY},\DDelta_{XY})$, and on this only. Both components $\Ss_{XY}$ and $\DDelta_{XY}$ may characterise a deviation from the state of independence, irrespective of the other: if $\Ss_{XY}$ is not rectangular, then $X$ and $Y$ cannot be independent even when $\DDelta_{XY} = \zero$ (Section \ref{subsec:regio}), and $\DDelta_{XY} \neq \zero$ precludes independence, even on a rectangular support (Corollary \ref{cor:indep}). Hence a proper quantification of the overall dependence must involve two contributions: one for the effect of the shape of the support $\Ss_{XY}$ (`{regional dependence}'), the other one assessing through $\DDelta_{XY}$ the deviation from `quasi-independence' on that support. When $\Ss_{XY}$ is rectangular, the first contribution should be 0 and the second would then measure deviation from independence as such. 

\ppn Below we explore this idea for the case of a bivariate discrete distribution on a finite support -- a more thorough treatment is left for future research. So we assume here that $\piXY$ is the distribution of a vector supported on $\Ss_{XY} \subseteq \Ss_{X} \times \Ss_Y = \{0,1,\ldots,R-1\} \times \{0,1,\ldots,S-1\}$ ($2 \leq R,S < \infty$). The distribution $\piXY$ may thus be identified to a non-negative matrix of $\Ms_{R \times S}$ with entries $p_{xy}$ summing to 1. 

\subsection{Regional dependence}

If $\Ss_{XY} = \Ss_{X} \times \Ss_Y$ (rectangular support), the dependence in $\piXY$ is entirely characterised by $(R-1)(S-1)$ parameters (the components of $\DDelta_{XY}$); see Section \ref{subsec:discr}$(a)$. By contrast, if the support is {not} rectangular (Section \ref{subsec:discr}$(b)$), then the dimension of $\DDelta_{XY}$ is $\dim(\Gamma_{\Ss_{XY}}^\circ)< (R-1)(S-1)$, where $\Gamma_{\Ss_{XY}}^\circ$ is the subspace of $\Ms_{R \times S}$ consisting of matrices with rows and columns summing to 0 and support contained in $\Ss_{XY}$. It is not so much that the dependence can be described by less parameters, it is rather that the support restrictions consume some of the initial $(R-1)(S-1)$ degrees of freedom. For example, once a support like (\ref{eqn:suppex}) or (\ref{eqn:suppdiag}) is identified for a $(3 \times 3)$-distribution, only one parameter suffices for completing the full description of dependence. Thus, compared to the $(R-1)(S-1) = 4$ necessary parameters on a rectangular support of that size, the supports (\ref{eqn:suppex}) and (\ref{eqn:suppdiag}) accounts for 3 dependence parameters. If the support is (\ref{eqn:suppdiag2}), all 4 initial dependence parameters are consumed by the support constraints. 

\ppn Distributions with rectangular support do not show any regional dependence, while for distributions such that $\DDelta_{XY}$ is empty (i.e., $\Gamma_{\Ss_{XY}}^\circ = \{\zero\}$), dependence is exclusively regional in nature. It seems, therefore, natural to quantify the amount of regional dependence in $\piXY$ by the fraction of dependence parameters consumed by the support constraints. Thus, we define
\[\textup{R}(X,Y) = \frac{(R-1)(S-1)-\dim(\Gamma_{\Ss_{XY}}^\circ)}{(R-1)(S-1)} \]	
as a measure of regional dependence. Evidently this index solely reflects (through $\Gamma_{\Ss_{XY}}^\circ$) the shape of $\Ss_{XY}$ inside $\Ss_X \times \Ss_Y$, but does not involve the distribution $\piXY$ as such -- in agreement with our description of regional dependence (Section \ref{subsec:regio}). Hence $\textup{R}(X,Y)$ is automatically invariant under the group transformation (\ref{eqn:equiv}), which preserves supports. Also, it is clear that $\textup{R}(X,Y) \in [0,1]$, with $\textup{R}(X,Y)=0 \iff \Ss_{XY} = \Ss_X \times \Ss_Y$ and $\textup{R}(X,Y)=1 \iff \Gamma_{\Ss_{XY}}^\circ = \{\zero\}$.

\begin{remark} \cite{Holland86} proposed a measure of regional dependence which, in the notation of this paper, would amount to $1 - (\piX \times \piY) (\Ss_{XY})$. However, this measure explicitly depends on the marginal distributions $\piX$ and $\piY$, and fails the basic test grounded in Corollary \ref{cor:excl}.	\qed
\end{remark}

\subsection{Deviation from quasi-independence on $\Ss_{XY}$}

Define $\Lambda_{\scriptscriptstyle XY} \in \Ms_{R \times S}$ as the matrix with component $(x,y) \in \{0,1,\ldots,R-1\} \times \{0,1,\ldots,S-1\}$ equal to\footnote{Here we arbitrarily set $\left[\Lambda_{\scriptscriptstyle XY}\right]_{xy} = 0$ for $(x,y) \not\in \Ss_{XY}$, but this is inconsequential.}
\[ \left[\Lambda_{\scriptscriptstyle XY}\right]_{xy} = \left\{\begin{array}{l l}  \log p_{xy} & \text{ if } (x,y) \in \Ss_{XY} \\ 0 & \text{ if } (x,y) \not\in \Ss_{XY}  \end{array}\right. .\]
Denote $d_\circ \doteq  \dim(\Gamma_{\Ss_{XY}}^\circ)$, which we assume positive here (the dependence is not exclusively regional), and consider an arbitrary  {orthonormal} basis $\{E_{k}: k = 1, \ldots,d_\circ \}$ of $\Gamma_{\Ss_{XY}}^\circ$. Then we can write $\DDelta_{XY} = (\delta_{{\scriptscriptstyle XY};1}, \ldots,\delta_{{\scriptscriptstyle XY};d_\circ})$ with $\delta_{{\scriptscriptstyle XY};k} \doteq \Lambda_{\scriptscriptstyle XY} \cdot E_k$ ($k=1,\ldots,d_\circ$) and $A \cdot B = \sum_{x,y} A_{xy} B_{xy}$ the Frobenius inner product between two matrices $A,B \in \Ms_{R \times S}$. Clearly, $\DDelta_{XY}$ is the vector of coordinates of the orthogonal projection, say $\Lambda^\circ_{\scriptscriptstyle XY}$, of $\Lambda_{\scriptscriptstyle XY}$ onto $\Gamma_{\Ss_{XY}}^\circ$; that is $\Lambda^\circ_{\scriptscriptstyle XY} = \sum_{k=1}^{d_\circ} \delta_{{\scriptscriptstyle XY};k} E_k$. If $X$ and $Y$ are quasi-independent on $\Ss_{XY}$, then $\DDelta_{XY} \equiv 0$ (Corollary \ref{cor:indep}), meaning that the projection of $\Lambda^\circ_{\scriptscriptstyle XY}$ is the null matrix. The Frobenius norm  
\begin{equation} \|\Lambda^\circ_{\scriptscriptstyle XY}\|_\text{F} = \sqrt{\sum_{k=1}^{d_\circ} \delta^2_{{\scriptscriptstyle XY};k} } = \|\DDelta_{XY}\|_2 \label{eqn:Frob} \end{equation}
can therefore be thought of as quantifying how remote from quasi-independence is the dependence structure of $\piXY$ on $\Ss_{XY}$. The norm $\|\DDelta_{XY}\|_2$ is invariant to a change of (orthonormal) basis. Note that the bases $\{E^{(00)}_{xy}\}$ and $\{E^{(\text{loc})}_{xy}\}$ leading to the usual sets of odds ratios (Section \ref{subsec:discr}) are {not} orthonormal, but orthonormal versions may easily be obtained via Gram-Schmidt. Other `interpretable' bases may be constructed as in \cite{Egozcue15} and \cite{Favcevicova18}.

\begin{example} \label{ex:bivBernmeas} Consider a bivariate Bernoulli distribution with full support ($p_{xy}>0\ \forall (x,y)$) as in Example \ref{ex:bivBern}. Here $d_\circ = 1$ with (normalised) basis element $\frac{1}{2} \left\{e_{00} + e_{11} - e_{10} - e_{01}\right\}$, in the notation of Section \ref{subsec:discr}. Hence $\DDelta_{XY} = \log \frac{\sqrt{p_{00}p_{11}}}{\sqrt{p_{10}p_{01}}} = \frac{1}{2} \log \omega$, and $\|\DDelta_{XY}\|_2 = \frac{1}{2} |\log \omega|$. For ease of interpretation, it may be desirable to transform monotonically this value to $[0,1]$. E.g., we could take $\tanh(\|\DDelta_{XY}\|_2/2) = |\Upsilon| \in [0,1]$, where $\Upsilon$ is Yule's colligation coefficient \citep{Yule12} (see \citet[Section 5.5]{Geenens2020} for reasons why this choice may be natural). The dependence is properly {quantified} by an absolute value, but {characterising} it entirely requires the signed number $\Upsilon$: the distribution cannot be fully identified from $p_{\bullet 1}, p_{1 \bullet}$ and $|\Upsilon|$ only. \qed
\end{example}

\noindent When $d_\circ > 1$, $\|\DDelta_{XY}\|_2$ provides a single-number quantification of the magnitude of the multiple entries which characterise the part of the dependence of $\piXY$ which is not regional. We may want to consider an increasing transformation $T: \R^+ \to [0,1)$, with $T(0) = 0$, and propose 
\begin{equation} \textup{Q}(X,Y) = T(\|\DDelta_{XY}\|_2) \in [0,1) \label{eqn:Q} \end{equation}
as a measure of deviation from quasi-independence on $\Ss_{XY}$. As a function of $\DDelta_{XY}$ only, $\textup{Q}(X,Y)$ is invariant under the group transformation (\ref{eqn:equiv}), and clearly, $\textup{Q}(X,Y) = 0 \iff$ $X$ and $Y$ are quasi-independent. The transformation $T$ is arbitrary and only aims at calibration.

\subsection{An overall measure of dependence}

Combining the previous two contributions, one may suggest the following overall measure of dependence: 
\begin{equation} \textup{D}(X,Y) = \textup{R}(X,Y) + (1-\textup{R}(X,Y)) \textup{Q}(X,Y). \label{eqn:mesdep} \end{equation}
This is the weighted average of $\textup{Q}(X,Y)$ and 1 with weights $1-\textup{R}(X,Y)$ and $\textup{R}(X,Y)$, respectively. In a sense, this amounts to setting the `missing' entries of $\DDelta_{XY}$ (those dependence parameters consumed by the support constraints) to their maximal value. 

\ppn This may be justified by referring to Example \ref{ex:bivBernmeas}. As all probabilities $p_{xy}$ are positive, the support is rectangular so that $\textup{R}(X,Y) = 0$ and $\textup{D}(X,Y) = \textup{Q}(X,Y) = |\Upsilon|$ (using the above suggested transformation $T(\cdot) = \atanh(\cdot/2)$). If one of the probabilities $p_{xy}$ approaches 0 (while staying positive), then $\log \omega \to \pm \infty$ and $|\Upsilon| \to 1$, so that $\textup{D}(X,Y) \to 1$. Now, suppose that one of the probabilities $p_{xy}$ is effectively 0: then $\Gamma_{XY}^\circ = \{\zero\}$, and $\textup{R}(X,Y) = 1$. The induced dependence is entirely regional, and consequently $\textup{D}(X,Y) = 1$, too. The proposed measure (\ref{eqn:mesdep}) is, therefore, continuous in the probabilities $p_{xy}$, including at 0.

\ppn Further properties of $\textup{D}(X,Y)$ easily follow from those of $\textup{R}(X,Y)$ and $\textup{Q}(X,Y)$. In particular, $\textup{D}(X,Y)$ is invariant under the group transformation (\ref{eqn:equiv}), as both $\textup{R}(X,Y)$ and $\textup{Q}(X,Y)$ are, and hence properly describes {dependence}: $\textup{D}(X,Y)$ is constant over any class $[\!\![\dep ]\!\!]$. Also, $\textup{D}(X,Y)$ fulfils the relevant requirements of \cite{Renyi59}, qualifying valid dependence measures:

\begin{itemize}
	\setlength{\itemsep}{1pt}
	\setlength{\parskip}{0pt}
	\setlength{\parsep}{0pt}
	\item[(A)] {$\textup{D}(X,Y)$ is defined for any pair of random variables $X$ and $Y$, neither of them being constant with probability 1};
	\item[(B)] $\textup{D}(X,Y) = \textup{D}(Y,X)$;
	\item[(C)] $0 \leq \textup{D}(X,Y) \leq 1$;
	\item[(D)] $\textup{D}(X,Y) = 0 \iff X \indep Y$;
	\item[(E)] $\textup{D}(X,Y) =1$ {if there is a strict dependence between $X$ and $Y$, i.e., either $X= f(Y)$ or $Y = g(X)$, where $f$ and $g$ are Borel-measurable functions};
	\item[(F)] {If the Borel-measurable functions $f$ and $g$ map the real axis in a one-to-one way onto itself, $\textup{D}(f(X),g(Y))$ $= \textup{D}(X,Y)$};
	%\item If the joint distribution of $X$ and $Y$ is normal, then $\textup{D}(X,Y) = |\rho|$, where $\rho$ is Pearson's correlation between $X$ and $Y$.
\end{itemize}

\noindent Properties (A),\footnote{Admittedly we only consider discrete variables, here.} (B) and (C) are obvious. Property (D) follows as $\textup{D}(X,Y) = 0 \iff \textup{R}(X,Y) = 0$ and $\textup{Q}(X,Y) = 0 \iff \Ss_{XY} = \Ss_X \times \Ss_Y$ and $\DDelta_{XY} =0 \iff X \indep Y$, by Corollary \ref{cor:indep}. What \cite{Renyi59} called `strict dependence' in (E) here reduces to the situation where there is exactly one positive probability $p_{xy}$ for each $x \in \Ss_X$ if $R \geq S$ or exactly one positive probability $p_{xy}$ for each $y \in \Ss_Y$ if $S \geq R$ (if $R=S$ then there is a one-to-one relationship between $X$ and $Y$). In that case $\Gamma_{XY}^\circ = \{\zero\}$, hence $\textup{R}(X,Y) = 1 $ and $\textup{D}(X,Y) = 1$. Property (F) was mentioned in Section \ref{sec:concord}, and in terms of $\textup{D}(X,Y)$ is reflected by the fact that $\dim(\Gamma_{\Ss_{XY}}^\circ)$ and hence $\textup{R}(X,Y)$ are invariant under permutations of the row and/or columns of $\piXY$, and $\textup{Q}(X,Y)$ is invariant to changes of the (orthonormal) basis $\{\gamma_1,\gamma_2,\ldots\}$ (which also covers permutations of the row and/or columns of $\piXY$). Note that \cite{Renyi59}'s last requirement is:
\begin{itemize}\item[(G)] {If the joint distribution of $X$ and $Y$ is normal, then $\textup{D}(X,Y) = |\rho|$, where $\rho$ is Pearson's correlation between $X$ and $Y$}.
\end{itemize}
Evidently this does not apply to the above discrete setting. This said, a continuous analogue of $\textup{Q}(X,Y) = \textup{D}(X,Y)$ (as the support of the bivariate Gaussian is rectangular, we would have $\textup{R}(X,Y) =0$) may easily be devised. Like in Section \ref{subsec:cont}$(a)$, consider the representation of a bivariate Gaussian distribution with correlation $\rho \in (-1,1)$ as the Gaussian copula density $c_\rho$ on $[0,1]^2$.  With the `amputated' Haar basis considered in Section \ref{subsec:cont} (which happens to be orthonormal in $L_2([0,1]^2)$), the entries of $\DDelta_{UV}$ are $\delta_{{\scriptscriptstyle UV};k {\bm \ell}} = \iint \log c_\rho \, \psi_{k {\bm \ell}}  \diff u \diff v$, for $k \in \N_0$ and ${\bm \ell} \in \{0,\ldots,2^k-1\}^2$. Hence the continuous analogue to (\ref{eqn:Frob}) would be
\[\|\DDelta_{UV}\|^2_2 = \iint_{[0,1]^2} \bar{\lambda}^2_{\rho}\, \diff u \diff v,   \]
where, similarly to (\ref{eqn:lambdaXY}), we have 
\[\bar{\lambda}_{\rho}  \doteq \log c_\rho - \int \log c_\rho \, \diff u - \int \log c_\rho \, \diff v + \iint \log c_\rho \, \diff u \diff v, \]
the centred log-odds ratio function of the copula density $c_\rho$. As $\log c_\rho \in L_2([0,1]^2)$ for all $\rho \in (-1,1)$, $\|\DDelta_{UV}\|_2 < \infty$, and can be verified to be an increasing function of $|\rho|$, indeed. With the appropriate transformation $T$ in (\ref{eqn:Q}), we may make it happen that $\textup{D}(U,V)=\textup{Q}(U,V) = |\rho|$, fulfilling \cite{Renyi59}'s postulate (G). This is similar to some ideas recently explored in \cite{Genest22}.

\section{Conclusion} \label{sec:ccl}

Close to 50 years ago, \citet[p.\,262]{Mosteller77} already regretted the `{troublesome, if not dangerous, confusion}' surrounding the usage of the word `dependence' in statistics. The comment remains very topical, as `dependence' and its declension continue to be used confusedly for referring to both the state of `non-independence' of two variables {\it and} the quantifiable strength of the mutual influence between them necessarily induced by that state. If the former unambiguously describes the negative of the precisely defined `independence' concept, a proper definition of the latter seems to be lacking. As coherently attacking many important problems requires the dependence between variables to be plainly understood as a specific structure of which we should be able to provide a precise and nuanced description, the lack of any clear framework around this leads to undesirable looseness and subjectivity when addressing essential questions. In response, this paper proposes a general definition for the dependence between two random variables defined on the same probability space.   

\ppn If one is willing to accept that definition (Definition \ref{def:dep}), then a cascade of implications establishes that the dependence $\dep$ in a bivariate distribution $\piXY$ admits a general representation (Corollary \ref{cor:dep}), `universal' in the sense that it is valid irrespective of the nature of the variables $X$ and $Y$ at play. This representation consists of two elements. The first one accounts for `regional dependence', the contribution to $\dep$ due to a non-rectangular joint support for $\piXY$, reflecting the incompatibility between some values of $X$ and some values of $Y$. The second one expresses in which way specific values of $X$ are more or less likely to occur with specific values of $Y$ on that support, and can be generally described as an `odds-ratio-like' object (or equivalent). The universal representation allows unequivocal identification of the proper way of capturing the dependence in particular situations of interest. Known results are recovered in familiar cases, but novel perspectives arise for less well-studied ones.

\ppn The proposed framework opens up numerous avenues for future research. For example, here we have considered exclusively bivariate distributions, that is, we focused on the dependence between two scalar variables. Though, Definition \ref{def:dep} can be extended naturally to the very general situation of  a collection of random objects for which we may define a joint distribution. This would allow the identification of dependence between an arbitrary number of random vectors, for example. Also, Section \ref{sec:depmes} explored the idea of a dependence measure based on the obtained representation -- thus being consistent with the given definition of dependence -- in the simplest case of a bivariate discrete distribution with finite support. An important question is how to extend this to make it more `universal'. 

%\ppn Finally, it seems essential to elucidate the nature of the connection between the ideas developed here and copulas. According to Sklar's theorem \citep{Sklar59}, every bivariate distribution $\piXY$ may be understood as a {coupling} $C$ of its margins $\piX$ and $\piY$. This coupling, characterised by a {copula}, has been largely equated to the dependence between $X$ and $Y$. Yet, if obvious parallels may be drawn between some of the observations made here and the rationale behind the copula approach, from the perspective of Sections \ref{sec:dep} and \ref{sec:universal0} the correspondence `dependence $\dep \sim $ copula $C$' is far from obvious. That connection is the topic of a follow-up paper \citep{Geenens23}.

\section*{Acknowledgments}

The author thanks Pierre Lafaye de Micheaux (UNSW Sydney, Australia) and Ivan Kojadinovic (Universit\'e de Pau et des Pays de l'Adour, France) for useful comments and discussions.

\bibliographystyle{../../elsarticle-harv-cond-mod}
\setlength{\bibsep}{0cm}
\def\bibfont{\footnotesize} %if class is 11 \footnotesize = 9pt $\small = 10pt http://en.wikibooks.org/wiki/LaTeX/Formatting#Sizing_text
\bibliography{../../libraries-copula}

\appendix 

\section{Definition of the set $\overline{\powerset}_{\Xss \times \Yss}$} \label{app:A}

Let $\powerset_{\Xss \times \Yss}$ be the set of all probability distributions $\piXY$ on $\Xss \times \Yss$ such that $\piXY \ll \mu_\Xss \times \mu_\Yss$. Let $\phiXY = \diff \piXY / \diff (\mu_\Xss \times \mu_\Yss)$ be the corresponding density, $\Ss_{XY}$, $\Ss_X$ and $\Ss_Y$ be the joint and marginal supports, as described in Section \ref{subsec:fw}. The subset $\overline{\powerset}_{\Xss \times \Yss}$ consists of the distributions $\piXY \in \powerset_{\Xss \times \Yss}$ such that:
\begin{enumerate}[(1)]
	\item $\phiXY \in L_\infty(\Xss\times \Yss,\mu_\Xss\times \mu_\Yss)$;
	\item \label{req:2} there exist rectangles $\{(E_i,F_i) \in \Bs_\Xss \otimes \Bs_\Yss\}_{i=1}^m$, $\{(G_j,H_j) \in \Bs_\Xss \otimes \Bs_\Yss\}_{j=1}^n$ of positive $(\mu_\Xss \times \mu_\Yss)$-measure, such that $\Ss_X = \bigcup_{i=1}^m E_i$ and $\Ss_Y = \bigcup_{j=1}^n G_j$ (up to sets of measure 0); and a constant $c>0$ such that $\phiXY(x,y) \geq c$ $(\mu_\Xss \times \mu_\Yss)\text{-a.e.}$ over $\bigcup_{i=1}^m (E_i \times F_i)$ and $\bigcup_{j=1}^n (G_j \times H_j)$;
	\item \label{req:3} there exists a sequence of measurable sets $K_j \in \Bs_\Xss \otimes \Bs_\Yss$ ($j=1,2,\ldots$) and a sequence of scalars $\varepsilon_j >0$ such that $\lim_{j\to\infty} \varepsilon_j =0$, with 
	\begin{enumerate}[(a)]
		\item $\phiXY >\varepsilon_j$ $(\mu_\Xss \times \mu_\Yss)\text{-a.e.}$ over $K_j$;
		\item $K_j$ is a countable union of measurable rectangles;
		\item $\Ss_{XY} = \bigcup_{j=1}^\infty K_j$.
	\end{enumerate}
\end{enumerate}
These basic requirements are reasonably mild. For example, they are automatically fulfilled if $\Ss_X$ and $\Ss_Y$ are compact and $\phiXY$ is continuous \citep[Theorem 2.34]{Nussbaum93} -- this is the case of any discrete distribution with finite support, among others. In our context, the above constraints may even be thought of as being much milder, in the following sense: we are actually allowed flexibility for choosing the working density $\phiXY$. Indeed we are free to pick any convenient `canonical' marginal supports and any one-to-one mappings $\Phi$ and $\Psi$ for moving to those supports (Section \ref{subsec:equiv}), producing any hand-picked `nice' marginal distributions $\piX$ and $\piY$ which may be convenient. In addition, the measures $\mu_\Xss$ and $\mu_\Yss$ (Section \ref{subsec:fw}) are left free, too; so that the density itself may be taken with respect to any convenient reference. In fact, Corollary 2.31 and Remark 3.14 in \cite{Nussbaum93} guarantee that, if Lemma \ref{lem:nussbaum} is valid for a distribution $\piXY$, then it would also hold true for any distribution $\piXYstar$ such that $\piXYstar \approx \piXY$ and $\phiXYstar \in L_\infty(\Xss\times \Yss,\mu_\Xss\times \mu_\Yss)$. So, all what is required is to find a joint measure equivalent to (some marginally-transformed version of) $\piXY$ with a density with respect to `some' product measure satisfying the requirements. This flexibility is justified by the fact that the above choices do not disturb the dependence of a given distribution.

\section{Proofs of the main results} \label{app:proofs}

\subsection*{Proof of Lemma \ref{lem:Leuridan}}

%(The proof is inspired by that of \cite{Brossard18}'s Theorem 3.1.)

Any distribution $\piXYstar$ in $\powerset_{\Xss \times \Yss} \cap \Sigma_{\piXY} \cap \Fs(\piXstar,\piYstar)$ would have a density $\phiXYstar = \diff \piXYstar/ \diff (\mu_\Xss \times \mu_\Yss)$ such that
\begin{equation} \left\{\begin{array}{l} \int \phiXYstar \diff \mu_\Yss = \phiXstar \qquad \text{($\piX$-a.e.)}\\ 
		\int \phiXYstar \diff \mu_\Xss = \phiYstar \qquad \text{($\piY$-a.e.)} \\ 
		\iint \phiXYstar \indic{\Ss^c_{XY}} \diff (\mu_\Xss \times \mu_\Yss) = 0 \end{array} \right. , \label{eqn:MOM} \end{equation}
where $\phiXstar = \diff \piXstar / \diff \mu_\Xss$, $\phiXstar = \diff \piYstar / \diff \mu_\Yss$ and $\Ss^c_{XY}$ is the complement in $\Xss \times \Yss$ of the support $\Ss_{XY}$ of $\piXY$. Define the closed convex cone $C = \{ \varphi \in L_1(\Xss \times \Yss,\mu_\Xss \times \mu_\Yss): \varphi \geq 0 \text{ a.e.}\} \subset L_1(\Xss \times \Yss,\mu_\Xss \times \mu_\Yss)$, and let $\As$ be the linear operator from $\Xi \doteq L_1(\Xss \times \Yss,\mu_\Xss \times \mu_\Yss)$ to $\Upsilon \doteq (L_1(\Xss,\mu_\Xss),L_1(\Xss,\mu_\Xss),L_1(\Yss,\mu_\Yss),L_1(\Yss,\mu_\Yss),\R)$ defined as 
\[ \As \varphi = \left(\int \varphi \diff \mu_\Yss,-\int \varphi \diff \mu_\Yss,\int \varphi \diff \mu_\Xss,-\int \varphi \diff \mu_\Xss,\iint\varphi \indic{\Ss^c_{XY}}\diff (\mu_\Xss \times \mu_\Yss)  \right).\]
The existence of a solution to (\ref{eqn:MOM}) means that $\bb \doteq (\phiXstar,-\phiXstar,\phiYstar,-\phiYstar,0) \in \Upsilon$ belongs to the (closed convex) cone $\As(C) \doteq \left\{\upsilon \in \Upsilon: \exists \varphi \in C \text{ s.t. } \upsilon = \As\varphi \right\}$. Following Farkas lemma (see, e.g., \citet[Theorem 2.2.6]{Craven78}), we have the equivalence 
\begin{equation} \bb \in \As(C) \iff \left(u \in \Upsilon', \As^T u \in C' \Longrightarrow u(\bb) \geq 0\right), \label{eqn:Farkas} \end{equation}
where $\Upsilon'$ is the dual of $\Upsilon$, $\As^T$ is the adjoint operator of $\As$ and $C'$ is the dual cone of $C$, i.e., $C' = \{\xi \in \Xi': \xi(\varphi) \geq 0 \ \forall\, \varphi \in C\}$, with $\Xi'$ the dual of $\Xi$.

\ppn Now, $\As^T u \in C' \iff (\As^T u)(\varphi) \geq 0$ for all $\varphi \in C$ $\iff$  $u(\As\varphi) \geq 0$ for all $\varphi \in C$, where 
\begin{align*} u(\As\varphi)\quad  =\quad  & \iint a^{(+)} \varphi \diff \mu_\Yss \diff \mu_\Xss - \iint a^{(-)} \varphi \diff \mu_\Yss \diff \mu_\Xss \\ & + \iint b^{(+)} \varphi \diff (\mu_\Xss \times \mu_\Yss)  - \iint b^{(-)} \varphi \diff (\mu_\Xss \times \mu_\Yss)  \\ & + k \iint_{\Ss_{XY}^c} \varphi \diff (\mu_\Xss \times \mu_\Yss) \end{align*}
for some $a^{(+)} ,a^{(-)}  \in L_\infty(\Xss,\mu_\Xss)$, $b^{(+)} ,b^{(-)}  \in L_\infty(\Yss,\mu_\Yss)$ and $k \geq 0$. Setting $a = a^{(+)} -a^{(-)}$, $b=b^{(+)} -b^{(-)}$, we have
\[ u(\As\varphi) = \iint (a + b + k \indic{\Ss_{XY}^c})\, \varphi \diff (\mu_\Xss \times \mu_\Yss).\]
Likewise, 
\[ u(\bb) = \int a\phiXstar \diff \mu_\Xss + \int b \phiYstar \diff \mu_\Yss. \]

\ppn $(a)$ $\powerset_{\Xss \times \Yss} \cap  \Sigma_{\piXY} \cap \Fs(\piXstar,\piYstar) = \emptyset$ means that there is no solution to (\ref{eqn:MOM}), that is, $\bb \not\in \As(C)$, which by (\ref{eqn:Farkas}) is equivalent to the existence of some $u \in \Upsilon'$ such that $\As^T u \in C'$ and $u(\bb) <0$. From above, 
\begin{align} \As^T u \in C' & \iff  \iint (a + b + k\indic{\Ss_{XY}^c})\, \varphi \diff (\mu_\Xss \times \mu_\Yss) \geq 0 \ \forall\,\varphi \in C \notag \\ & \iff a(x)+b(y) +k\indic{(x,y) \not\in \Ss_{XY}} \geq 0 \quad  (\mu_\Xss \times \mu_\Yss)\text{-a.e.} \label{eqn:abk} \end{align}
while 
\begin{equation} u(\bb) <0 \iff  \int a\phiXstar \diff \mu_\Xss + \int b \phiYstar \diff \mu_\Yss< 0. \label{eqn:ubneg}\end{equation}

\ppn Let $X^*$ (resp.\ $Y^*$) be a random variable defined on $\Xss$ (resp.\ $\Yss$) with distribution $\piXstar$ (resp.\ $\piYstar$), and define $U = a(X^*)$ and $V=-b(Y^*)$, for $a,b$ two functions satisfying (\ref{eqn:abk}) and (\ref{eqn:ubneg}). We have:
\[\int_{-\infty}^\infty (\P(U>t) -\P(V>t)) \diff t = \E(U) - \E(V) = \int a \diff \piXstar + \int b \diff \piYstar < 0 \]
(from (\ref{eqn:ubneg})). Thus $\mathfrak{T} \doteq \{t \in \R:\P(U>t) -\P(V>t)<0 \} \neq \emptyset$. For any $t \in \mathfrak{T}$, define $A_t \doteq \{x \in \Xss: a(x) \leq t\}$ and $B_t \doteq \{y \in \Yss: b(y) <-t\}$ so that for any $(x,y) \in A_t \times B_t$, $a(x) + b(y) <0$. From (\ref{eqn:abk}), such $(x,y)$ cannot be in $\Ss_{XY}$. Thus for any $t \in \mathfrak{T}$, $(A_t, B_t) \in \text{N}_{\piXY}$, while
\[\piXstar(A_t) + \piYstar(B_t) = \int_{A_t} \phiXstar \diff \mu_\Xss + \int_{B_t} \phiYstar \diff \mu_\Yss = \P(U\leq t)  + \P(V>t) = 1 +\P(V>t) - \P(U>t)>1.  \]

\ppn $(b)$ $(i)$ The claim is trivial, by the convexity of the sets $\Sigma_\piXY$ and $\Fs(\piXstar,\piYstar)$.

\ppn $(ii)$  ($\Leftarrow$) Let $(A,B) \in \text{N}_{\piXY}$ such that $\piXstar(A) + \piYstar(B) = 1$ and $(A^c,B^c) \not\in \text{N}_{\piXY}$. For any $\piXYstar \in \powerset_{\Xss \times \Yss} \cap \Sigma_{\piXY} \cap \Fs(\piXstar,\piYstar)$ (including $\overline{\pi}^*_{\scriptscriptstyle XY}$ from $(b)(i)$), $\piXYstar(A \times B^c)=\piXstar(A)$ and $\piXYstar(A^c \times B^c)= \piYstar(B^c) - \piXYstar(A \times B^c) = \piXstar(A)-\piXYstar(A \times B^c) = 0$. Hence, as $\piXY(A^c \times B^c) \neq 0$ by assumption, $\piXY \not\ll \piXYstar$ and they are not measure-theoretically equivalent. 

\ppn $(ii)$  ($\Rightarrow$) By definition we have $\overline{\pi}^*_{\scriptscriptstyle XY} \ll \piXY$, but assume that $ \piXY \not\ll \overline{\pi}^*_{\scriptscriptstyle XY}$. Then there exists a measurable set $S_0 \subset \Ss_{XY}$ such that $\overline{\pi}^*_{\scriptscriptstyle XY}(S_0) = 0$ while $\piXY(S_0) > 0$. By $(i)$, all distributions $\piXYstar \in \powerset_{\Xss \times \Yss} \cap   \Sigma_{\piXY} \cap \Fs(\piXstar,\piYstar)$ are such that $\piXYstar \ll \overline{\pi}^*_{\scriptscriptstyle XY}$, so $\piXYstar \in \powerset_{\Xss \times \Yss} \cap   \Sigma_{\piXY} \cap \Fs(\piXstar,\piYstar) \Rightarrow \piXYstar(S_0) = 0$. Now, $\piXYstar \in \powerset_{\Xss \times \Yss} \cap  \Sigma_{\piXY} \cap \Fs(\piXstar,\piYstar)$ means that its density $\phiXYstar$ satisfies (\ref{eqn:MOM}), of which $\iint \phiXYstar \indic{S_0} \diff (\mu_\Xss \times \mu_\Yss) = 0$ should therefore be a consequence. Thus the system (\ref{eqn:MOM}) augmented by the strict inequality 
\begin{equation} \iint \phiXYstar \indic{S_0} \diff (\mu_\Xss \times \mu_\Yss) > 0 \label{eqn:strict} \end{equation} must be inconsistent. As above, we can resort to Farkas lemma to show that such inconsistency is equivalent to the existence of some functions $a \in L_\infty(\Xss,\mu_\Xss)$, $b  \in L_\infty(\Yss,\mu_\Yss)$, $k \geq 0$ and $\ell >0$ such that 
\begin{equation} a+ b+ k \indic{\Ss_{XY}^c} - \ell \indic{S_0} = 0 \qquad (\mu_\Xss \times \mu_\Yss)\text{-a.e.} \label{eqn:abkl} \end{equation}
and 
\begin{equation} \int a \diff \piXstar + \int b \diff \piYstar \leq 0. \label{eqn:ubneg2} \end{equation}
(Note that the inequality is not strict here, as opposed to (\ref{eqn:ubneg}), as a consequence of the strict inequality in (\ref{eqn:strict}). Also, $\ell >0$; {\it \`a la} \citet[Theorem 4.2.4]{Webster94}.) 

\ppn Let $X^*$ (resp.\ $Y^*$) be a random variable defined on $\Xss$ (resp.\ $\Yss$) with distribution $\piXstar$ (resp.\ $\piYstar$), and define $U = a(X^*)$ and $V=-b(Y^*)$, for $a,b$ two functions satisfying (\ref{eqn:abkl}) and (\ref{eqn:ubneg2}). We have:
\[\int_{-\infty}^\infty (\P(U>t) -\P(V>t)) \diff t = \E(U) - \E(V) = \int a \diff \piXstar + \int b \diff \piYstar \leq 0. \]
Denote $u_0$ the lower bound of the support of $U$, and $v_0$ the upper bound of the support of $V$ (these bounds are finite as $a,b$ are essentially bounded). As $\P(U>t) -\P(V>t) \geq 0$ for $t < u_0$ or $t \geq v_0$, there necessarily exists $t \in [u_0,v_0)$ such that $P(U>t) -\P(V>t) \leq 0$. For such $t$, define $A_t \doteq \{x \in \Xss: a(x) \leq t\}$ and $B_t \doteq \{y \in \Yss: b(y) <-t\}$ so that for any $(x,y) \in A_t \times B_t$, $a(x) + b(y) <0$ (a.e.), implying that 
\[ k \indic{\Ss_{XY}^c} - \ell \indic{S_0} > 0\]
from (\ref{eqn:abkl}). As $S_0 \subset \Ss_{XY}$, it follows that $(x,y)$ cannot belong to $\Ss_{XY}$ when $(x,y) \in A_t \times B_t$; that is, $\piXY(A_t \times B_t) = 0$, and $(A_t,B_t) \in \text{N}_{\piXY}$. In addition,
\[\piXstar(A_t) = \P(U \leq t) = \P(U \in [u_0,t)) > 0 \quad \text{ and } \quad \piYstar(B_t) = \P(V > t) = \P(V \in (t,v_0]) > 0  \]
and 
\[\piXstar(A_t)+ \piYstar(B_t) =  \P(U \leq t)+ \P(V > t) = 1 - \P(U>t) +\P(V > t) \geq 1. \]
It must be that $\piXstar(A_t)+ \piYstar(B_t)  = 1$, as by assumption $\piXstar(A) + \piYstar(B) \leq 1$ $\forall (A,B) \in \text{N}_{\piXY}$. 

\ppn As $(A_t \times B_t) \subset S_{XY}^c$ and $S_0 \subset S_{XY}$, $(A_t \times B_t) \subset S_0^c$, that is, $S_0 \subset (A_t \times B_t)^c = (A_t \times B_t^c) \cup (A_t^c \times B_t) \cup (A_t^c \times B_t^c)$. Assume that $(A_t^c,B_t^c) \in \text{N}_{\piXY}$. Then $\Ss_{XY} \subseteq (A_t \times B_t^c) \cup (A_t^c \times B_t)$. As $\piXstar(A_t)+ \piYstar(B_t)  = 1$, we know that $\piXstar(A_t)= \piYstar(B^c_t)$, and $\piXYstar(A_t \times B_t^c) = \piXstar(A_t)$. Then there cannot be $\tilde{A} \subset A_t$ and $\tilde{B} \subset B^c_t$ such that $(\tilde{A},\tilde{B}) \in \text{N}_{\piXY}$ and $\piXstar(\tilde{A}) + \piYstar(\tilde{B}) > \piXstar(A_t)$, otherwise $(\tilde{A},\tilde{B} \cup B_t) \in \text{N}_{\piXY}$ and $\piXstar(\tilde{A})+ \piYstar(\tilde{B}) + \piYstar(B_t) > \piXstar(A_t) + \piYstar(B_t) = 1$, which is not allowed by assumption. So the restriction of $\piXY$, $\piXstar$ and $\piYstar$ may be handled in the same way as the above argument. If for all $\tilde{A} \subset A_t$ and $\tilde{B} \subset B^c_t$ such that $(\tilde{A},\tilde{B}) \in \text{N}_{\piXY}$, $\piXstar(\tilde{A}) + \piYstar(\tilde{B}) < \piXstar(A_t)$, then $\overline{\pi}^*_{\scriptscriptstyle XY} \approx \piXY$ on $A_t \times B_t^c$ and thus $S_0 \not\subset A_t \times B_t^c$ -- we come to the same conclusion recursively if there exists $\tilde{A} \subset A_t$ and $\tilde{B} \subset B^c_t$ such that $(\tilde{A},\tilde{B}) \in \text{N}_{\piXY}$, $\piXstar(\tilde{A}) + \piYstar(\tilde{B}) = \piXstar(A_t)$. Similarly, one shows that $S_0 \not\subset A_t^c \times B_t$. But if $(A_t^c,B_t^c) \in \text{N}_{\piXY}$, $\piXY(A_t^c  \times B_t^c) = 0$ and $S_0$ cannot be in $A_t^c  \times B_t^c$, leading to a contradiction. Hence $(A_t^c,B_t^c) \not\in \text{N}_{\piXY}$.

\ppn Thus $ \piXY \not\ll \overline{\pi}^*_{\scriptscriptstyle XY}$ implies the existence of a set $(A,B) \in \text{N}_{\piXY}$ such that $\piXstar(A) + \piYstar(B) = 1$, and $(A^c,B^c) \not\in \text{N}_{\piXY}$.

\subsection*{Proof of Proposition \ref{prop:compa}}

($\Rightarrow$) Assume that $\{\phiYgX\}$ and $\{\phiXgYstar\}$ are compatible. By definition, there exists a bivariate distribution $\widetilde{\pi}_{\scriptscriptstyle XY}$ admitting $\{\phiYgX\}$ and $\{\phiXgYstar\}$ as conditional densities. Thus $\{\phiYgX\}$ and $\{\phiXgYstar\}$ must encapsulate the same dependence structure, which is the dependence structure of $\widetilde{\pi}_{\scriptscriptstyle XY}$. Hence $\Delta(\{\phiYgX\}) = \Delta(\{\phiXgYstar\})$.

\ppn ($\Leftarrow$) Assume that $\Delta(\{\phiYgX\}) = \Delta(\{\phiXgYstar\})$. By (\ref{eqn:depsup}), this implies that any distribution admitting $\{\phiYgX\}$ as $Y|X$-conditional density shares the same joint support (and hence marginal supports) as any other distribution admitting $\{\phiXgYstar\}$ as $X|Y$-conditional density. In particular, take a distribution $\piX$ on $\Ss_X$ such that $\diff \piX/\diff \mu_\Xss \in L_\infty(\Xss,\mu_\Xss)$, and a distribution $\piYstar$ on $\Ss_Y$ such that $\diff \piYstar/\diff \mu_\Yss \in L_\infty(\Yss,\mu_\Yss)$, and form 
\[ 	\piXY(\cdot)  \doteq \iint \phiYgX \indic{\cdot} \diff (\piX \times \mu_\Yss) \quad \text{ and } \quad  \piXYstar(\cdot)  \doteq \iint \phiXgYstar \indic{\cdot} \diff \piYstar \diff \mu_\Xss. \]
As these two distributions have conditional densities $\{\phiYgX\}$ and $\{\phiXgYstar\}$, respectively, they have the same support. See from Corollary \ref{cor:Leuridan2} that condition (\ref{eqn:compacond}) is automatically fulfilled for $\piXY$ and the marginals $\piXstar, \piYstar$ of $\piXYstar$, given that $\textup{N}_{\piXY} = \textup{N}_{\piXYstar}$ and there exists a distribution in $\Fs(\piXstar,\piYstar)$ with the suitable support. Then, according to Corollary \ref{cor:alphabetadep}, there is a unique distribution in $\Fs(\piXstar,\piYstar)$ with the same dependence as $\piXY$ -- this must be $\piXYstar$, by assumption -- and it is such that 
\[\frac{\diff \piXYstar}{\diff \piXY} = \alpha^*_{\scriptscriptstyle X} \beta^*_{\scriptscriptstyle Y}\]
for some positive functions $\alpha^*_{\scriptscriptstyle X} \in L_1(\Xss,\mu_\Xss)$ and $\beta^*_{\scriptscriptstyle Y} \in L_1(\Yss,\mu_\Yss)$. This implies that 
\[\frac{ \phiXgYstar}{ \phiYgX} = \frac{\diff \piXYstar}{\diff (\mu_\Xss \times \piYstar)} \frac{\diff (\piX \times \mu_\Yss)}{\diff \piXY} = \left(\alpha^*_{\scriptscriptstyle X}\frac{\diff \piX}{\diff \mu_\Xss}\right) \left(\beta^*_{\scriptscriptstyle Y}\frac{\diff \mu_\Yss}{\diff \piYstar}\right) \doteq \widetilde{\alpha}^*_{\scriptscriptstyle X}\widetilde{\beta}^*_{\scriptscriptstyle Y} . \]
As $\diff \piX/\diff \mu_\Xss \in L_\infty(\Xss,\mu_\Xss)$, $\widetilde{\alpha}^*_{\scriptscriptstyle X} \in L_1(\Xss,\mu_\Xss)$. \citet[Theorem 4.1]{Arnold89} establishes this as a necessary and sufficient condition for $\{\phiXgYstar\}$ and $\{\phiYgX\}$ to be compatible. \qed

\subsection*{Proof of Theorem \ref{thm:equivclass}}

Let $\{\phiYgX\}$ be the $Y|X$-conditional density of $\piXY$, and $\{\phiXgYstar\}$ be the $X|Y$-conditional density of $\piXYstar$. As per Proposition \ref{prop:compa}, $\dep = \dep^*$ if and only if $\{\phiYgX\}$ and $\{\phiXgYstar\}$ are compatible. From \citet[Theorem 4.1]{Arnold89}, $\{\phiYgX\}$ and $\{\phiXgYstar\}$ are compatible if and only if $(i)$  $\{(x,y) \in \Ss_X  \times \Yss:   \phiYgX(y|x) > 0\} = \{(x,y) \in \Xss  \times \Ss_Y:   \phiXgYstar(x|y) > 0\}$, which is equivalent to $\Ss_{XY} = \Ss_{XY}^*$ and therefore follows from (\ref{eqn:depsup}); and $(ii)$ there exist (a.e.-)positive functions $\phi_{\scriptscriptstyle X} \in L_1(\Xss,\mu_\Xss)$ and $\phi_{\scriptscriptstyle Y} \in L_1(\Yss,\mu_\Yss)$, such that
\begin{equation*} \frac{\phiXgYstar}{\phiYgX} = \frac{\phi_{\scriptscriptstyle X}}{\phi_{\scriptscriptstyle Y}} \qquad \text{a.e.\ on $\Ss_{XY}$}. \end{equation*} 
As 
\[\frac{\phiXgYstar}{\phiYgX} = \frac{\diff \piXYstar}{\diff (\mu_\Xss \times \pi^*_Y)} \frac{\diff (\piX \times \mu_\Yss)}{\diff \piXY} \]
directly from (\ref{eqn:conddens}), we have
\begin{equation} \frac{\diff \piXYstar}{\diff \piXY} =  \frac{\phi_{\scriptscriptstyle X} \diff \mu_\Xss}{\diff \piX} \frac{\diff \piYstar}{\phi_{\scriptscriptstyle Y} \diff \mu_\Yss} \doteq \alpha^*_{\scriptscriptstyle X} \beta^*_{\scriptscriptstyle Y}. \label{eqn:dpidpii}\end{equation}
See that $\alpha^*_{\scriptscriptstyle X} = \phi_{\scriptscriptstyle X} \diff \mu_\Xss / \diff \piX \in L_1(\Xss,\piX)$, as $\int \alpha^*_{\scriptscriptstyle X} \diff \piX = \int  \phi_{\scriptscriptstyle X} \diff \mu_\Xss  < \infty$. Likewise, $1/\beta^*_{\scriptscriptstyle Y} = \phi_{\scriptscriptstyle Y} \diff \mu_\Yss / \diff \piYstar \in L_1(\Yss,\pi^*_{\scriptscriptstyle Y})$. From (\ref{eqn:dpidpii}) we have directly that $\diff \piXY/\diff \piXYstar =  1/\alpha^*_{\scriptscriptstyle X} \times 1/\beta^*_{\scriptscriptstyle Y}$, while, swapping $\piXYstar$ and $\piXY$ from the beginning of the argument, we would have obtained instead of (\ref{eqn:dpidpii}): 
\[ \frac{\diff \piXY}{\diff \piXYstar} =  \frac{\phi^*_{\scriptscriptstyle X} \diff \mu_\Xss}{\diff \piXstar} \frac{\diff \piY}{\phi^*_{\scriptscriptstyle Y} \diff \mu_\Yss}\]
for some functions $\phi^*_{\scriptscriptstyle X} \in L_1(\Xss,\mu_\Xss)$ and $\phi^*_{\scriptscriptstyle Y} \in L_1(\Yss,\mu_\Yss)$. Identifying both expressions, we see that $\beta^*_{\scriptscriptstyle Y}$ can be written $\phi^*_{\scriptscriptstyle Y} \diff \mu_\Yss / \diff \piY$, showing that $\beta^*_{\scriptscriptstyle Y} \in L_1(\Yss,\piY)$. The claim follows. % from that we have shown:
%\[  \{\phiXgY\} \text{ and } \{\phiYgXstar\} \text{ are compatible } \iff (i) \ \Ss_{XY} = \Ss^*_{XY} \text{ and } (ii)\ \frac{\diff \piXYstar}{\diff \piXY} =  \alpha_X \beta_Y, \]
%for a $\piX$-integrable function $\alpha_X$ and a $\piY$-integrable function $\beta_Y$. 
\qed

\subsection*{Proof of Theorem \ref{thm:maxinv}}

We show that, for two distributions $\piXY,\piXYstar \in \overline{\powerset}_{\Xss \times \Yss}$,
\[\piXY \sim \piXYstar \iff \iint \log \phiXY \, \diff \gamma_k = \iint  \log \phiXYstar \, \diff \gamma_k \quad \forall k \in \{1,2,\ldots\}, \]
where $\piXY \sim \piXYstar$ is defined as per (\ref{eqn:equiv2}) and $\{\gamma_1,\gamma_2,\ldots\}$ is any basis of $\Gamma_{{XY}}^\circ$ defined in (\ref{eqn:Gamma}).

\ppn ($\Rightarrow$) Assume $\piXY \sim \piXYstar$. This means that $\phiXYstar = \phiXY \alpha^*_{\scriptscriptstyle X} \beta^*_{\scriptscriptstyle Y}$ for two (a.e.-)positive functions $\alpha^*_{\scriptscriptstyle X} \in L_1(\Xss,\piX)$ and $\beta^*_{\scriptscriptstyle Y} \in L_1(\Yss,\piY)$. Note that, by symmetry (see Proof of Theorem \ref{thm:equivclass}), $\alpha^*_{\scriptscriptstyle X}$ and $\beta^*_{\scriptscriptstyle Y}$ are also essentially bounded. Then, $\forall k \in \{1,2,\ldots\}$,
\[ \iint  \log \phiXYstar \, \diff \gamma_k  = \iint  \log \phiXY \, \diff \gamma_k + \iint  \log \alpha^*_{\scriptscriptstyle X} \, \diff \gamma_k + \iint  \log \beta^*_{\scriptscriptstyle Y} \, \diff \gamma_k  = \iint  \log \phiXY \, \diff \gamma_k \]
by definition, as any $\gamma_k \in \Gamma_{{XY}}^\circ$ has null marginals, making $\iint  \log \alpha^*_{\scriptscriptstyle X} \, \diff \gamma_k = \iint  \log \beta^*_{\scriptscriptstyle Y} \, \diff \gamma_k =0$.

\ppn ($\Leftarrow$) Assume that $\iint \log \phiXY \, \diff \gamma_k = \iint  \log \phiXYstar \, \diff \gamma_k$, $\forall k \in \{1,2,\ldots\}$. As $\{\gamma_1,\gamma_2,\ldots\}$ is a basis of $\Gamma_{{XY}}^\circ$, it follows that 
\[\iint \log \frac{\phiXYstar}{\phiXY}\, \diff \gamma = 0 \qquad \forall \gamma \in \Gamma_{{XY}}^\circ.\]
Hence $\log \frac{\phiXYstar}{\phiXY}$ must belong to the annihilator \citep[Section 4.6]{Rudin91} of $\Gamma_{{XY}}^\circ$, viz.
$$\Psi_{XY}^\circ \doteq \{\psi \in L_\infty(\Xss \times \Yss): \iint \psi \, \diff \gamma = 0, \ \forall \gamma \in \Gamma_{{XY}}^\circ\},$$
which can be verified to be $\Psi_{XY}^\circ = L_\infty(\Xss,\mu_\Xss) \oplus L_\infty(\Yss,\mu_\Yss)$. Hence there exist $\psi_{\scriptscriptstyle X} \in L_\infty(\Xss,\mu_\Xss)$ and $\psi_{\scriptscriptstyle Y} \in L_\infty(\Yss,\mu_\Yss)$ such that $\log \frac{\phiXYstar}{\phiXY} = \psi_{\scriptscriptstyle X} + \psi_{\scriptscriptstyle Y}$, and finally $\frac{\diff \piXYstar}{\diff \piXY} = \frac{\phiXYstar}{\phiXY} = \exp(\psi_{\scriptscriptstyle X} + \psi_{\scriptscriptstyle Y}) \doteq \alpha^*_{\scriptscriptstyle X}\beta^*_{\scriptscriptstyle Y}$, where $\alpha^*_{\scriptscriptstyle X}\beta^*_{\scriptscriptstyle Y}>0$. As $\alpha^*_{\scriptscriptstyle X}$ and $\beta^*_{\scriptscriptstyle Y}$ are essentially bounded, they are $\piX$- and $\piY$-integrable. \qed

\end{document}